\documentclass[a4paper]{iopart}
\pdfoutput=1
\usepackage{etex}
\expandafter\let\csname equation*\endcsname\relax 
\expandafter\let\csname endequation*\endcsname\relax 
\usepackage{amssymb}
\usepackage{amsmath}
\usepackage{graphicx}
\usepackage{subfig}

\usepackage[inline]{asymptote}

\makeatletter
\newsavebox{\sfe@box}
\newenvironment{subfloatenv}[1]{%
\def\sfe@caption{#1}%
\setbox\sfe@box\hbox\bgroup\color@setgroup}%
{\color@endgroup\egroup\subfloat[\sfe@caption]%
{\usebox{\sfe@box}}}
\makeatother

\newcommand{\norm}[1]{\ensuremath{\lVert#1\rVert}}
\newcommand{\grade}[1]{\ensuremath{\langle#1\rangle}}

\newcommand{\reverse}[1]{\ensuremath{\widetilde{#1}}}

\newcommand{\R}[1]{\ensuremath{\mathbb{R}^{#1}}}
\newcommand{\T}[1]{\ensuremath{\mathbb{T}^{#1}}}
\newcommand{\E}[1]{\ensuremath{\mathbb{E}^{#1}}}
\newcommand{\M}[1]{\ensuremath{{M}_{#1}}}

\newcommand{\El}[1]{\ensuremath{{E}_{#1}}}
\newcommand{\Hy}[1]{\ensuremath{{H}_{#1}}}

\newcommand{\dS}[1]{\ensuremath{{dS}_{#1}}}
\newcommand{\AdS}[1]{\ensuremath{{AdS}_{#1}}}

\renewcommand{\S}[1]{\ensuremath{\mathbb{S}^{#1}}}

\newcommand{\tb}[1]{\ensuremath{\textbf{#1}}}

\newcommand{\ts}[1]{\ensuremath{\mathsf{#1}}}

\renewcommand{\e}{\tb{e}} 
\newcommand{\I}{\tb{I}} 

\newcommand{\Id}{\ensuremath{\mathcal{R}}} 
\newcommand{\J}{\ensuremath{\mathcal{D}}}
\newcommand{\m}{\makebox[0pt][r]{\(-\)}}

\begin{document}
\title{A key to the projective model of homogeneous metric spaces}
\author{Andrey Sokolov}
\address{School of Physics, University of Melbourne, Parkville, VIC 3010, Australia}
\ead{a.sokolov@pgrad.unimelb.edu.au}

\begin{abstract}
A metric introduced on a projective space yields a homogeneous metric space known as a Cayley-Klein geometry.
This construction is applicable not only to Euclidean and non-Euclidean spaces but also to kinematic spaces (space-times).
A convenient algebraic framework for Cayley-Klein geometries called the projective model is 
developed in \cite{gunn2011homogeneous,gunn2011geometry}.
It is based on Grassmann and Clifford algebras and provides a set of algebraic tools for
modeling points, lines, planes and their geometric transformations such as projections and isometries.
Isometry groups and their Lie algebras find a natural and intuitive expression in the projective model.
The aim of this paper is to translate the foundational concepts of the projective model
from the language of projective geometry to a more familiar language of vector algebra and thereby
facilitate its spread and adoption among physicists and applied mathematicians.
I apply the projective model to Minkowski, de-Sitter, and anti de-Sitter space-times in two dimensions.
In particular, I show how the action of the Poincar\'e group can be captured by the Clifford algebra
in a uniform fashion with respect to rotations (boosts) and translations.
\end{abstract}

\submitto{}



\section{Introduction}
Prior to Descartes' introduction of coordinates, Euclidean geometry was treated from a synthetic perspective, according to which
points, lines, and planes are considered given and their geometric relations are governed by a set of axioms \cite{stillwell2005four}.
Descartes' analytic geometry bridged a gap between geometry and algebra by enabling one to apply algebraic techniques to solve geometric problems.
The concept of a vector as a single entity described by coordinates emerged gradually during the nineteenth century.
Although matrix notation for linear transformations was introduced by Cayley in the late 1850s,
modern definition of vector spaces and linear transformations appeared only towards the end of the nineteenth century.
Quaternions were invented by Hamilton in 1843 and later proposed as a preferred geometric language for the three-dimensional space,
which is important in physical sciences.
The product of two quaternions whose scalar component is zero can be reduced to two operations: the dot and cross product of three-dimensional vectors.
Based on this decomposition, Gibbs devised vector algebra and analysis, which became popular towards the end of the nineteenth century
and spread to physical sciences through their use in electromagnetic theory popularised by Heaviside
(for detailed history see \cite{crowe2002history}).
Vectors eventually became the language of choice in most fields of classical physics including Newtonian mechanics, electrodynamics, and fluid dynamics.
Since vector algebra and analysis are limited to three-dimensional space, 
they are supplemented by a more general tensor algebra when dealing with four-dimensional space-time in special and general relativity.

Vectors in \R{n} can be visualised as oriented line segments stemming from the origin.
As such they are suitable for representing lines passing through the origin of \R{n}, which can be identified with one-dimensional linear subspaces of \R{n}.
This basic picture can be extended in two different ways with the help of Grassmann algebra and projective geometry.
Grassmann algebra was first introduced in 1844 and in revised form in 1862, but its value had not been realised until the early twentieth century
in the works of Cartan.
It extends the concept of vectors by introducing new algebraic entities such as bivectors and trivectors for representing general linear subspaces of \R{n}
 (see \cite{dorst2009geometric} for an elementary introduction).
This algebra is used widely in differential geometry, where it is known as exterior algebra
with bivectors and trivectors appearing in the guise of differential forms.
The extension provided by projective geometry (and closely related affine geometry) 
harkens back to the synthetic perspective on geometry, characteristic of classical Euclidean geometry,
as it provides techniques for modeling points, lines, and planes located anywhere in space
by means of the so-called homogeneous coordinates \cite{richtergebert2011perspectives}.
Projective geometry is different from affine geometry in that it also allows one to model points, lines, and planes at infinity,
which simplifies theorems by eliminating special cases and has other advantages.

Projective geometry is favoured in pure mathematics and serves as a foundation of algebraic geometry.
It is also used extensively in computer graphics \cite{shreiner2004opengl,vince2010mathematics}, 
where  homogeneous coordinates are employed to represent points and lines in three-dimensional space.
Nevertheless, it has rarely been used in physics and applied mathematics 
due to its abstract formulation and apparent lack of a convenient algebraic framework.
Moreover, projective geometry is non-metric and therefore requires an additional structure for incorporating metric properties.
This structure is provided by the so-called Cayley-Klein construction \cite{richtergebert2011perspectives,klein1939elementary,onishchik2006projective}, 
which was described first by Cayley in 1859 and developed in detail in early 1870s 
by Klein as part of his Erlangen program.
The choice of a metric in the Cayley-Klein construction decides the type of metric geometry one obtains from projective geometry.
This gives rise to a range of homogeneous metric spaces called Cayley-Klein geometries, which serve as models of Euclidean and non-Euclidean spaces,
including popular kinematic spaces such as Minkowski space-time.
The Cayley-Klein construction is quite elaborate both in its set-up and implementation.
Without an efficient and intuitive framework for computations, Cayley-Klein geometries will remain unpopular.
Ironically, the algebra for just such a framework was invented by Clifford at about the same time as Klein's work on Cayley-Klein geometries,
but the connection has not been recognised until recently \cite{gunn2011homogeneous,gunn2011geometry}.

Clifford algebra emerges as a natural extension of Grassmann algebra that incorporates metric properties of \R{n} 
\cite{lounesto2001clifford,garling2011clifford,hestenes1987clifford,meinrenken2013clifford}.
It subsumes complex numbers and quaternions and was initially developed as a generalisation of these associative algebras.
Thanks to its connection with Grassmann algebra, it is closely linked to geometric structures of \R{n} 
and in fact can serve, like tensor algebra, as a generalisation of Gibbs'  vector algebra to spaces of arbitrary dimension.
Clifford algebra can be defined as the quotient of the tensor algebra and 
is therefore not as general as  tensor algebra, 
but due to its narrower scope it is more effective at representing geometric structures and their transformations.
It relies solely on exterior and orthogonal properties of linear subspaces of \R{n} and abstracts away all other properties carried by tensor algebra.
Clifford algebra is used in physics in the guise of Dirac matrices \cite{d2003geometric,francis2005construction}, 
differential geometry \cite{lawson1989spin,benn1987introduction}, Clifford analysis \cite{delanghe1992clifford,brackx1982clifford},
which is a generalisation of complex analysis to spaces of arbitrary dimension, and robotics \cite{selig2005geometric}.

Points, lines, and planes of the \(n\)-dimensional projective space can be identified with  linear subspaces of \R{n+1}
and therefore Grassmann algebra of \R{n+1} becomes the algebra of points, lines, and planes of the \(n\)-dimensional space.
The associated Clifford algebra inherits this structure from the Grassmann algebra.
In addition, it furnishes space with the metric properties,
consistent with the Cayley-Klein construction \cite{gunn2011geometry}, 
and allows one to represent the isometry group by elements of the algebra.
The group of isometries depends on the metric and consists of the  elements that act in various ways on points, lines, and planes, while preserving distances and angles.
The representation of isometries in the Clifford algebra compares favourably with the way isometries are represented by the classical matrix groups 
(see \cite{stillwell2008naive} for the background on matrix groups),
since elements of the matrix groups are generally difficult to interpret geometrically as they explicitly depend on coordinates,
whereas the representation in Clifford algebra is independent of coordinates.
The elements of the Clifford algebra that represent isometries preserving orientation constitute a Lie group,
whose Lie algebra can be represented by the same Clifford algebra in a geometrically transparent way.
The projective model is a combination of projective geometry and Clifford algebra that is consistent with the Cayley-Klein construction.
It embodies a synthesis of the old synthetic and modern analytic approaches to geometry.

The aim of this paper is to provide a concise introduction to the projective model developed in \cite{gunn2011homogeneous,gunn2011geometry},
without overburdening the reader with technical details that are not relevant for applications in physics and computing.
I develop the projective model for two-dimensional Cayley-Klein geometries only, 
as this allows for easy visualisation but still conveys most essential features of the model.
I give a brief introduction to the relevant aspects of projective geometry, 
such as projective duality and embedding into model vector spaces, in Section~\ref{section.projective}.
This section serves as the foundation for the following Section~\ref{section.grassmann} on the relevant Grassmann algebras,
where I describe how points and lines of the two-dimensional space are represented by elements of the Grassmann algebras
and briefly explore non-metric aspects of geometry in two dimensions.
In Section~\ref{section.clifford}, I introduce the Clifford algebra and explore metric aspects of geometry in two dimensions
for both kinematic (Minkowski, de-Sitter, anti de-Sitter) and non-kinematic spaces (Euclidean, hyperbolic, elliptic).

I omit the Cayley-Klein construction, since distances and angles can be defined within the Clifford algebra.
Instead, I concentrate on various geometric constructions and transformations provided by the Clifford algebra.
I reveal how space-time intervals can be obtained in the projective model without appeal to integration required in
the traditional approaches based on curved manifolds and pseudo-Riemannian geometry.
The Clifford algebra can be used to represent rotations in Minkowski space by conjugation with the elements of the algebra,
that is one gains access to the action of the Lorentz group on space-time events purely in terms of the Clifford algebra.
The projective model extends this application of the Clifford algebra to the Poincar\'e group as it can
represent not only rotations but also translations by conjugation.

\section{Projective foundations}
\label{section.projective}
The projective model relies on projective geometry.
The most relevant concepts for the following exposition are the projective duality and the related concept of the top-down model of geometry.
They are summarised and illustrated in this section for the case of the planar geometry (two-dimensional projective geometry).
A modern and comprehensive, yet accessible, treatment of projective geometry is given in \cite{richtergebert2011perspectives},
while \cite{klein1939elementary} provides a more elementary introduction.
I will use \T{2} to designate the target space whose geometry I wish to study, \T{2*} will designate 
the projective dual of the target space (or simply the dual space).
\R{3} will be used for the vector-space model of the target space (or simply the model space)
and \R{3*} for the usual vector-space dual of the model space.
Note that even though I employ identical notation for the projective dual
and the vector-space dual, projective duality is not equivalent to vector-space duality as will become apparent in the following.

\subsection{Projective duality}
Projective duality can be introduced in the following way (cf.\ \cite[p.\ 59]{klein1939elementary} and \cite[p.\ 57]{richtergebert2011perspectives}, 
see also \cite{conradt2000principle}).
Any line \(L\) in \T{2}, which does not pass through the origin, can be defined by
\begin{equation}
1+ax+by=0,
\label{line in T2}
\end{equation}
where \(a\) and \(b\) are some fixed real numbers and \(x\) and \(y\) are the coordinates of  points on the line.
Given Equation~(\ref{line in T2}), I can selects a particular line in \T{2} by specifying a pair of real numbers \((a,b)\).
The dual space \T{2*} can be thought of as the space of  pairs \((a,b)\) that define lines in \T{2} via Equation~(\ref{line in T2});
the meaning of the origin \((a,b)=(0,0)\) will be clarified below.
I will refer to the pairs such as \((a,b)\) as points of the dual space \T{2*}.
The point \((a,b)\in\T{2*}\) and the line \(L\subset\T{2}\) defined by~(\ref{line in T2}) are said to be dual to each other.
For instance, the line defined by \(1+3x+y=0\) is dual to the point \((3,1)\in\T{2*}\).

Equation~(\ref{line in T2}) is symmetric with respect to \((a,b)\) and \((x,y)\).
One can view \(x\) and \(y\) as fixed  and \(a\) and \(b\) as variable.
This defines a line \(K=\{(a,b)|1+ax+by=0\}\) in the dual space \T{2*}.
The line \(K\) is determined by specifying a point in \T{2} with the fixed coordinates \(x\) and \(y\).
Consequently, the target space \T{2} can be viewed as 
space of the pairs \((x,y)\) that define lines in \T{2*} via  Equation~(\ref{line in T2}).
The point \((x,y)\in\T{2}\) and the line \(K\subset\T{2*}\) are said to be dual to each other.
For instance, the line defined by \(1-a-\tfrac{1}{2}b=0\) is dual to the point \((-1,-\tfrac{1}{2})\in\T{2}\) (see Figure~\ref{line sheaf}).

\begin{figure}[t!]
\hspace{-0.3cm}%
\begin{subfloatenv}{\T{2}}
\begin{asy}
import Drawing2D;
Drawing2D drawing = Drawing2D(4, 0.1, unit_size=0.75);
drawing.grid();
drawing.target_axes();

pair n = normalise((-1,2));
pair c = (0.8,0.4);
labeled_pair[] points = points_on_line(center=c, direction=n, separation=1.5*sqrt(5), N=1);
labeled_pair[] extra_points = {labeled_pair((1,0),'(1,0)'), labeled_pair((0,2), '(0,2)')};
points.append(extra_points);

pair intersection = (-1,-1/2);
real[] relative_position = { 0.1, 0.1, 0.9, 0.1, 0.9};
align[] label_align = {NoAlign, NoAlign, (-0.2,1), NoAlign, NoAlign,};
for(int i=0; i<points.length; ++i) { 
  pair n = perp(points[i].pair);
  path l = line(center=intersection, direction=n, extent_from_center=4); 
  real th = degrees(n);
  if (i==0) { th -= 180; } 
  write(points[i].label, ' ', (string)th);
  draw(rotate(th)*Label('$\J'+points[i].label+'$', relative_position[i], align=label_align[i]), l); 
}
dot(intersection);
label('$\J(K)$', intersection, SE);
drawing.crop();
\end{asy}
\end{subfloatenv}
\begin{subfloatenv}{\T{2*}}
\begin{asy}
import Drawing2D;
Drawing2D drawing = Drawing2D(4, 0.1, unit_size=0.75);
drawing.grid();
drawing.dual_axes();

pair n = normalise((-1,2));
pair c = (0.8,0.4);
labeled_pair[] points = points_on_line(center=c, direction=n, separation=1.5*sqrt(5), N=1);
labeled_pair[] extra_points = {labeled_pair((1,0),'(1,0)'), labeled_pair((0,2), '(0,2)')};
points.append(extra_points);
align[] point_label_align = {W,W,W,SW+(0.1,0),W};
for(int i=0; i<points.length; ++i) { dot(Label('$'+points[i].label+'$', align=point_label_align[i]), points[i].pair);}
path l = line(center=c, direction=n, extent_from_center=5); 
draw(Label('$K$',0.6), l);

drawing.crop();
\end{asy}
\end{subfloatenv}
\caption{The line \(K=\{(a,b)|1-a-\tfrac{1}{2}b=0\}\) in the dual space (b) and the corresponding sheaf of lines in the target space (a).}
\label{line sheaf}
\end{figure}

Points are dual to lines and vice versa in two-dimensional projective geometry.
Hence, Equation~(\ref{line in T2}) defines the duality transformation denoted by \(\J\) acting on points and lines of the dual space \T{2*}
and yielding the corresponding lines and points of the target space \T{2}, e.g.\ \(\J(3,1)= \{(x,y)|1+3x+y=0\}\) and 
\(\J(\{(a,b)|1-a-\tfrac{1}{2}b=0\})=(-1,-\tfrac{1}{2})\).
I also define a transformation denoted by \(\Id\) that acts on points and lines of the dual space \T{2*}
and yields identical points and lines in the target space \T{2}, e.g.\ 
\(\Id(3,1)= (3,1)\) and 
\(\Id(\{(a,b)|1-a-\tfrac{1}{2}b=0\})=\{(x,y)|1-x-\tfrac{1}{2}y=0\}\).
The duality transformation \(\J\) turns points into lines and lines into points,
whereas the transformation \(\Id\) gives points for points and lines for lines.

To visualise the line \(L\) defined \((a,b)\) it is convenient to have an expression for at least one point on \(L\) in terms of the coordinates of the point \((a,b)\).
Fortunately, one such point can be readily obtained by intersecting \(L=\J(a,b)\) with 
the line that passes through the origin of \T{2} and the point \(\Id(a,b)\).
It will be called the central point of \(L\) and its coordinates are given by
\begin{equation}
x_c= \frac{-a}{a^2+b^2},\quad y_c = \frac{-b}{a^2+b^2}.
\end{equation}
For instance, the central point of the line defined by \((a,b)=(3,1)\) is located at \((-0.3,-0.1)\).
A similar definition can be given for the central point of a line in the dual space \(\T{2*}\).
Since only incidence of points and intersection of lines are used in the definition of the central point,
it is independent of the metric.

If the target space is Euclidean, the central point of the line \(L=\J(a,b)\) is a point on \(L\)  at the shortest distance from the origin,
and the distance from \(L\) to the origin is equal to the inverse of the distance
from the origin to \(\Id(a,b)\).
For instance, the distance from \(L=\J(3,1)\) to the origin is equal to \(\sqrt{(-0.3)^2+(-0.1)^2}=1/\sqrt{10}\)
and the distance from the origin to \(\Id(3,1)\) is equal to \(\sqrt{3^2+1^2}=\sqrt{10}\).
If the metric is not Euclidean, these relations are not applicable.

\subsubsection{Top-down model of geometry}
In the standard picture of geometry, points play the fundamental role and other geometric objects such as lines are built from points. 
For instance, it takes two points to define a line and the line is seen as a set of points which lie on the line.
I call this the bottom-up model of \T{2}, since higher-dimensional objects are constructed from points, which are zero-dimensional
and thus are located at the bottom of the hierarchy. 
This contrasts with the top-down model of \T{2} where lines, which are at the top of the hierarchy in a two-dimensional space,
are considered fundamental and 
lower-dimensional objects such as points are built from lines as their intersection.
Hence, a point is a derivative object defined by any two lines passing through the point.
Just as a line is seen as a set of points in the bottom-up model,
a point in the top-down model should be seen as a set of lines passing through the point.
I will refer to this set as a sheaf of lines attached to the point where the lines intersect.
The bottom-up and the top-down models can be applied to the dual space \T{2*} in the same fashion.

\begin{figure}[t!]
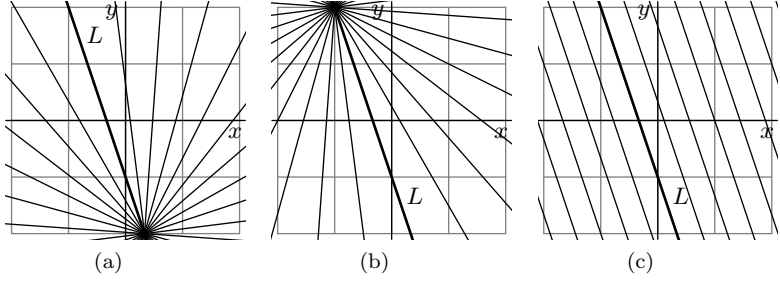

\hspace{1cm}%
\begin{subfloatenv}{ }
\begin{asy}
import Drawing2D;
Drawing2D drawing = Drawing2D(2, 0.1, unit_size=0.75);
drawing.grid();
drawing.target_axes();

pair n = (-1/3,1);
pair c = (-1/3,0);
path l = line(center=c, direction=n, extent_from_center=5); 
draw(Label("$L$", 0.65), l, currentpen+1);

pair[] ns = {};
int N=16;
real phi0 = angle(n);
for(int i=1; i<N; ++i) { pair n=expi(i*pi/N + phi0); ns.push(n); }

pair P = (1/3,-2);
for(pair n: ns) { path l = line(center=P, direction=n, extent_from_center=5);  draw(l); }

drawing.crop();
\end{asy}
\end{subfloatenv}%
\begin{subfloatenv}{ }
\begin{asy}
import Drawing2D;
Drawing2D drawing = Drawing2D(2, 0.1, unit_size=0.75);
drawing.grid();
drawing.target_axes();

pair n = (-1/3,1);
pair c = (-1/3,0);
path l = line(center=c, direction=n, extent_from_center=5); 
draw(Label("$L$", 0.35), l, currentpen+1);

pair[] ns = {};
int N=16;
real phi0 = angle(n);
for(int i=1; i<N; ++i) { pair n=expi(i*pi/N + phi0); ns.push(n); }

pair P = (-1,2);
for(pair n: ns) { path l = line(center=P, direction=n, extent_from_center=5);  draw(l); }

drawing.crop();
\end{asy}
\end{subfloatenv}%
\begin{subfloatenv}{ }
\begin{asy}
import Drawing2D;
Drawing2D drawing = Drawing2D(2, 0.1, unit_size=0.75);
drawing.grid();
drawing.target_axes();

pair n = (-1/3,1);
pair c = (-1/3,0);
path l = line(center=c, direction=n, extent_from_center=5); 
draw(Label("$L$", 0.35), l, currentpen+1);

pair[] ps = {};
int N=10;
for(int i=1; i<N; ++i) { pair P=c+(i,-i/3)/2; ps.push(P); }
for(int i=1; i<N; ++i) { pair P=c-(i,-i/3)/2; ps.push(P); }

for(pair P: ps) { path l = line(center=P, direction=n, extent_from_center=5);  draw(l); }

drawing.crop();
\end{asy}
\end{subfloatenv}
\caption{A line \(L\) in \T{2} and three points on \(L\) shown as sheaves in (a) and (b) and a stack in (c).}
\label{line and its point at infinity}
\end{figure}

The top-down model of geometry gives an instant access to an additional set of points which are not readily available in the bottom-up model.
Namely, besides points represented by sheaves of intersecting lines, one can also consider points represented by stacks of parallel lines as follows.
In the two-dimensional Euclidean space, every stack of parallel lines defines a specific point at infinity and
the lines that comprise the stack can be thought of as intersecting at that point.
This implies that a point at infinity can be approached by moving along any of the lines in the stack that defines the point.
Moreover, the same point at infinity is reached by moving in either direction along the lines in the stack.
In other words, from the point of view of projective geometry one observes the same point at infinity by looking in two opposite directions.
Any line \(L\) in the two-dimensional Euclidean space defines a stack of lines parallel to \(L\),
which can be interpreted as a point at infinity.
It is convenient to assume that this point at infinity lies \(L\) and therefore
the set of points comprising \(L\) is extended by the point at infinity (see Figure~\ref{line and its point at infinity}),
i.e.\ a line in projective space is topologically equivalent to a circle.
The set of all points at infinity may be called the line at infinity.

A stack of lines can be defined as a set of lines in which every line is defined by \(ax+by+d=0\),
where \(a\) and \(b\) are fixed and \(d\) spans all possible values in \R{}.
Since the definition of a stack of lines does not rely on any metric properties, it can be used not only 
in Euclidean space but in any target space \T{2} regardless of its metric if any.
A stack of lines in the target space \T{2} defines a point of the target space. 
However, this point might be at an infinite distance from the origin or the distance to it might even be undefined
in some metric spaces.
The terminology used in the previous paragraph is only applicable in Euclidean space,
but the points represented by stacks as well as the line that consists of such points are defined even if the target space is non-metric.
From this point of view, there is no difference between sheaves and stacks and, in fact,
both of these structures are usually referred to in projective geometry as pencils of lines.

The dual space  provides a convenient arena for the study of geometry of the target space  from the top-down point of view,
for one can apply a more familiar bottom-up model of the dual space \T{2*} in order to gain access to the top-down model of \T{2}.
I have demonstrated above that points of the dual space corresponds to lines of the target space via projective duality.
It turns out that lines of the dual space corresponds to sheaves or stacks of lines in the target space.
Indeed, according to the bottom-up model a line in the dual space is a set of points, each of which is dual to a certain line in the target space.
For example, consider a line \(K=\{(a,b)|1-a-\tfrac{1}{2}b=0\}\) shown in Figure~\ref{line sheaf}(b).
Every point \((a,b)\) on \(K\) corresponds to a line in the target space \T{2} via projective duality,
and it turns out that every such line passes through the point \(\J(K)\) (several representative lines are shown in Figure~\ref{line sheaf}(a)).
So, collectively points on \(K\) give rise to a sheaf of lines attached to the point \(\J(K)\), which is dual to \(K\).

\begin{figure}[t!]
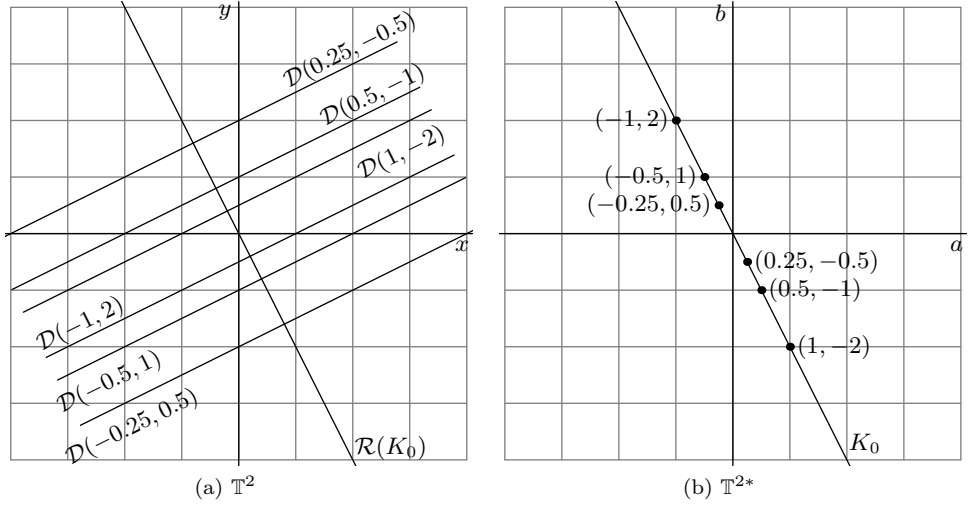

\hspace{-0.3cm}%
\begin{subfloatenv}{\T{2}}
\begin{asy}
import Drawing2D;
Drawing2D drawing = Drawing2D(4, 0.1, unit_size=0.75);
drawing.grid();
drawing.target_axes();

pair n = normalise((-1,2));
labeled_pair[] points = points_on_line(center=(0,0), direction=n, separation=sqrt(5)/2, N=2);

int index_of_origin = quotient(points.length-1, 2);
points.delete(index_of_origin); // delete (0,0)

labeled_pair[] extra_points = {labeled_pair((-0.25,0.5),'(-0.25,0.5)'), labeled_pair((0.25,-0.5), '(0.25,-0.5)')};
points.append(extra_points);

real[] relative_position = { 0.1, 0.1, 0.1, 0.1, 0.1, 0.1};
align[] label_align = {(0.5,-1), NoAlign, NoAlign, (-0.5,1), NoAlign, NoAlign};
for(int i=0; i<points.length; ++i) { 
  pair n = perp(points[i].pair);
  real a = points[i].pair.x;
  real b = points[i].pair.y;
  pair intersection = -1/5 * (1/a, 4/b);
  path l = line(center=intersection, direction=n, extent_from_center=4); 
  real th = degrees(n);
  if (th > 180) { th -= 180; }
  write(points[i].label, ' ', (string)th);
  draw(rotate(th)*Label('$\J'+points[i].label+'$', relative_position[i], align=label_align[i]), l); 
}

path l = line(center=(0,0), direction=n, extent_from_center=5); 
draw(Label('$\Id(K_0)$',0.07), l);

drawing.crop();
\end{asy}
\end{subfloatenv}
\begin{subfloatenv}{\T{2*}}
\begin{asy}
import Drawing2D;
Drawing2D drawing = Drawing2D(4, 0.1, unit_size=0.75);
drawing.grid();
drawing.dual_axes();

pair n = normalise((-1,2));
labeled_pair[] points = points_on_line(center=(0,0), direction=n, separation=sqrt(5)/2, N=2);

int index_of_origin = quotient(points.length-1, 2);
points.delete(index_of_origin); // betel (0,0)

labeled_pair[] extra_points = {labeled_pair((-0.25,0.5),'(-0.25,0.5)'), labeled_pair((0.25,-0.5), '(0.25,-0.5)')};
points.append(extra_points);

align[] point_label_align = {E,E,W,W,W,E};
for(int i=0; i<points.length; ++i) { dot(Label('$'+points[i].label+'$', align=point_label_align[i]), points[i].pair);}

path l = line(center=(0,0), direction=n, extent_from_center=5); 
draw(Label('$K_0$',0.07), l);

drawing.crop();
\end{asy}
\end{subfloatenv}
\caption{The line \(K_0=\{(a,b)|2a+b=0\}\) in the dual space (b) and the corresponding stack of lines in the target space (a).}
\label{stack of lines}
\end{figure}

Furthermore, it turns out that points on a line passing through the origin of \T{2*} give rise to a stack of lines in \T{2}.
Hence, in the two-dimensional Euclidean space, points at infinity are dual to lines passing through the origin of \T{2*},
and vice versa.
An example is shown in Figure~\ref{stack of lines} for \(K_0=\{(a,b)|2a+b=0\}\);
the point \(\J(K_0)\) is at infinity and can only be displayed as a stack of lines.
The origin of the dual space corresponds to the line at infinity in the target space
and therefore the line at infinity in \T{2} ought to be included in the definition of a stack of lines for consistency.
Observe also that in Euclidean space the line \(\Id(K_0)\) is perpendicular to any line in the stack collectively represented by points on \(K_0\).
In other metric spaces one can only claim that a line passing through the origin of \T{2*}
is dual to a point of \T{2} represented by a certain stack of lines, or more generally
any line in the dual space corresponds to a pencil of lines in the target space.

Identical constructions can be carried out for the top-down model of the dual space 
in terms of the bottom-up model of the target space. 
This in particular yields points in \T{2*} represented by stacks of lines in \T{2*}
and the line consisting of such points.
The former correspond to lines passing through the origin of \T{2} and the latter corresponds  to the origin of \T{2}.

\begin{figure}[t!]
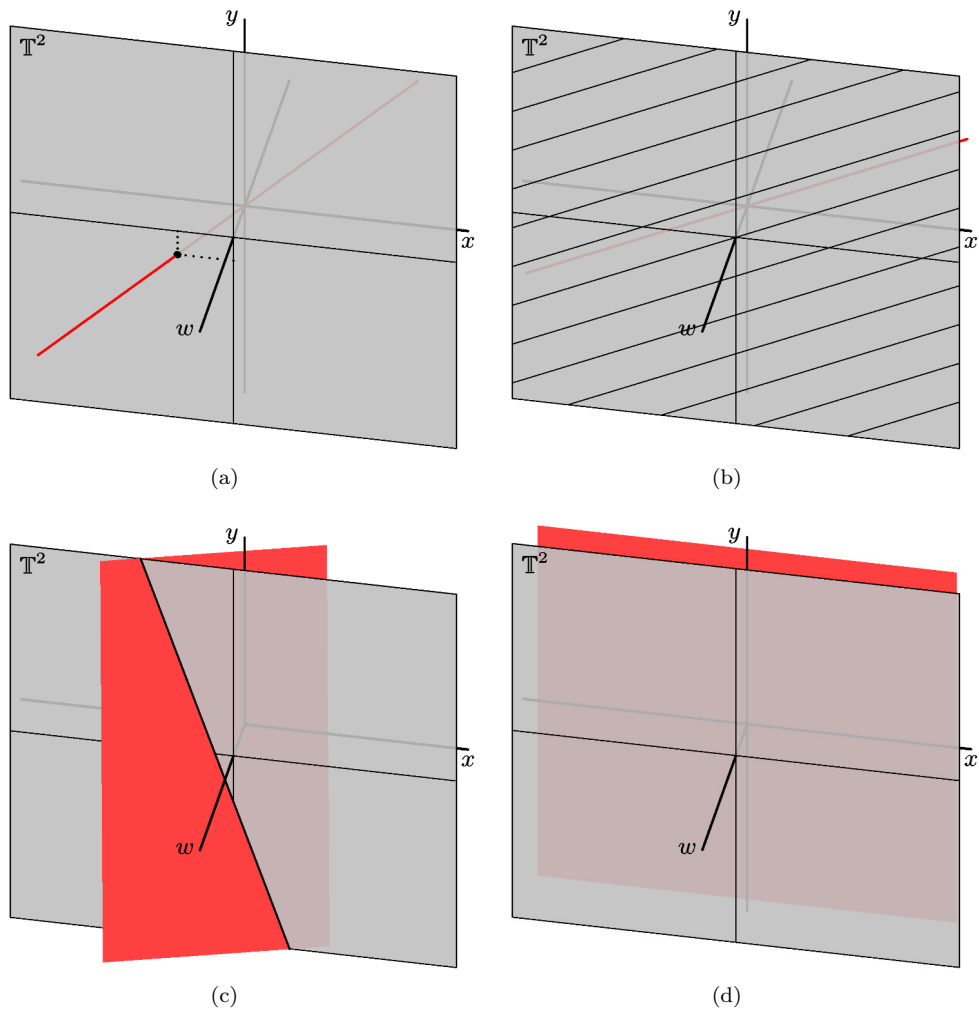

\begin{subfloatenv}{ }
\begin{asy}
import Drawing3D;
real e = 4;
DrawingR3 drawing = DrawingR3(e,0.1, axes_extent=e, camera=(10,2,7));
drawing.target_axes();

path3 p = plane_path3((1,0,0),(1,0,0), v=(0,1,0),scale=2*e);
drawing.Drawing3D.plane(p);
draw((1,-e,0)--(1,e,0));
draw((1,0,-e)--(1,0,e));

path3 l = line_path3((1,-1,-1/2)/4, direction=(1,-1,-1/2), extent_from_center=e+0.25);
draw(l,red+1);
dot((1,-1,-1/2));

draw((1,-1,0)--(1,-1,-1/2), Dotted);
draw((1,0,-1/2)--(1,-1,-1/2), Dotted);

label("$\mathbb{T}^{2}$", (1.5,-(e-0.5),e));
\end{asy}
\end{subfloatenv}
\begin{subfloatenv}{ }
\begin{asy}
import Drawing3D;
real e = 4;
DrawingR3 drawing = DrawingR3(e,0.1, axes_extent=e, camera=(10,2,7));
drawing.target_axes();

path3 p = plane_path3((1,0,0),(1,0,0), v=(0,1,0),scale=2*e);
drawing.Drawing3D.plane(p);
draw((1,-e,0)--(1,e,0));
draw((1,0,-e)--(1,0,e));

path3 l = line_path3((0,0,0), direction=(0,1,1/2), extent_from_center=e+0.4);
draw(l,red+1);

path3 l = line_path3((1,0,0), (0,1,1/2), 10);
path3 p = drawing.T2_path3(2*e);
int N=6;
for(int i=-N; i<=N; ++i) { path3 shifted_l = shift(-i*(0,-1,2)/3)*l; triple[] ps = intersectionpoints(shifted_l,p); draw(ps[0]--ps[1]); }

label("$\mathbb{T}^{2}$", (1.5,-(e-0.5),e));
\end{asy}
\end{subfloatenv}\hfill%
\begin{subfloatenv}{ }
\begin{asy}
import Drawing3D;
real e = 4;
DrawingR3 drawing = DrawingR3(e,0.1, axes_extent=e, camera=(10,2,7));
drawing.target_axes();

path3 p = plane_path3((1,0,0),(1,0,0), v=(0,1,0),scale=2*e);
drawing.Drawing3D.plane(p);
draw((1,-e,0)--(1,e,0));
draw((1,0,-e)--(1,0,e));

triple c = (1,-1/3,0);
path3 Q = plane_path3((1,3,1), c, theta=15, scale=2*e-0.9); 
drawing.Drawing3D.plane(Q,mediumred,nullpen,opacity=1);

path3 l = plane_intersection(Q,p); 
draw(l,black+1);

label("$\mathbb{T}^{2}$", (1.5,-(e-0.5),e));
\end{asy}
\end{subfloatenv}
\begin{subfloatenv}{ }
\begin{asy}
import Drawing3D;
real e = 4;
DrawingR3 drawing = DrawingR3(e,0.1, axes_extent=e, camera=(10,2,7));
drawing.target_axes();

path3 p = plane_path3((1,0,0),(1,0,0), v=(0,1,0),scale=2*e);
drawing.Drawing3D.plane(p);
draw((1,-e,0)--(1,e,0));
draw((1,0,-e)--(1,0,e));

triple c = (0,0,0);
path3 Q = plane_path3((1,0,0), c, theta=45, scale=2*e-0.5); 
drawing.Drawing3D.plane(Q,mediumred,nullpen,opacity=1);

label("$\mathbb{T}^{2}$", (1.5,-(e-0.5),e));
\end{asy}
\end{subfloatenv}
\caption{Embedding of the target space \T{2} into the model space \R{3}}
\label{embedding T2 in R3}
\end{figure}

\subsubsection{Embedding}
Unlike points of the target space represented by sheaves, points represented by stacks cannot be written in the form \((x,y)\).
This limitation is overcome in projective geometry by introducing the so-called homogeneous coordinates,
which is equivalent to 
identifying points of \T{2} with one-dimensional linear subspaces\footnote{Recall that a one-dimensional linear subspace is a line, 
which passes through the origin of \R{3}, and a two-dimensional linear subspace is a plane, 
which passes through the origin of \R{3}.} 
of \R{3} as follows.
The model space \R{3} consists of the triples \((w,x,y)\), where \((1,x,y)\in\R{3}\) is identified with \((x,y)\in\T{2}\),
so that a point \((x,y)\) of the target space can be identified with 
the one-dimensional linear subspace of \R{3} that passes through \((1,x,y)\).
An example is shown in Figure~\ref{embedding T2 in R3}(a) for \((x,y)=(-1,-\tfrac{1}{2})\), whose homogeneous coordianates are \((1,-1,\tfrac{1}{2})\) or any non-zero scalar multiple of that.
One can think of \T{2}, bar points represented by stacks, as a plane embedded into \R{3} at \(w=1\).
To define a point of the target space represented by a stack, one needs to specify any line which belongs to the stack,
e.g.\ \(ax+by=0\) (the other lines of the stack are then given by \(ax+by+d=0\), where \(d\) spans \R{}).
The equation \(ax+by=0\) defines a line in the target space and 
at the same time it defines a plane in the model space (by letting \(w\) range over \R{}),
whose intersection with the plane \(w=0\) is a one-dimensional linear subspace.
This subspace is confined to the plane \(w=0\) by construction.
It does not intersect the plane \(w=1\) at any point \((1,x,y)\) and it can thus be identified with the point in \T{2} represented by the stack.
An example is shown in Figure~\ref{embedding T2 in R3}(b), where the stack representing the point is defined by the line \(x-2y=0\) and its homogeneous coordinates are \((0,-2,1)\) or any non-zero scalar multiple of that.

The treatment of lines is similar to the treatment of points above.
Thanks to the embedding, a line \(L\) in the target space can be identified with a line in the plane \(w=1\).
It in turn can be identified with a two-dimensional linear subspace
of \R{3} that intersects the plane \(w=1\) along \(L\) 
(see Figure~\ref{embedding T2 in R3}(c), where \(L\) is defined by \(1+3x+y=0\)).
Furthermore, the line that consists of the points represented by stacks, e.g. the line at infinity in Euclidean space,
can be identified with the two-dimensional linear subspace \(w=0\) (see Figure~\ref{embedding T2 in R3}(d)).

\begin{figure}[t!]
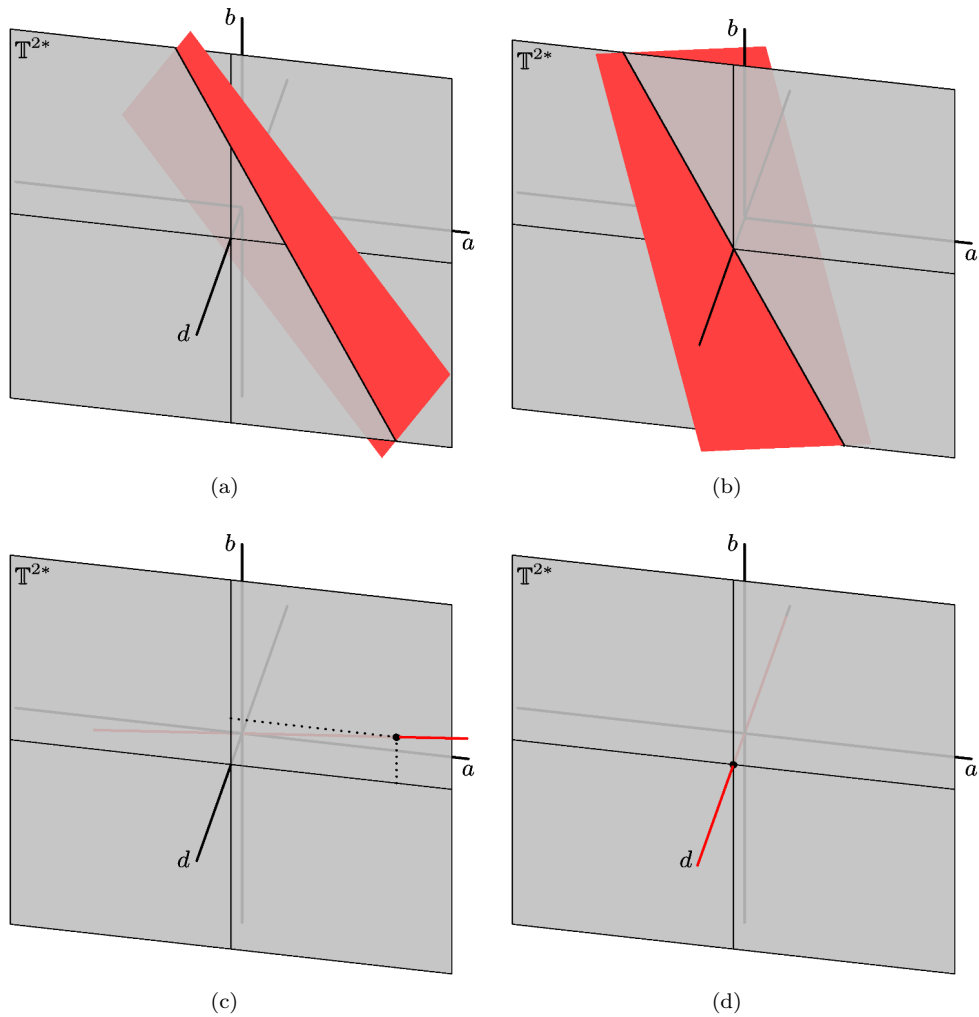

\begin{subfloatenv}{ }
\begin{asy}
import Drawing3D;
real e = 4;
DrawingR3 drawing = DrawingR3(e,0.1, axes_extent=e, camera=(10,2,7));
drawing.dual_axes();

path3 p = plane_path3((1,0,0),(1,0,0), v=(0,1,0), scale=2*e);
drawing.Drawing3D.plane(p);
draw((1,-e,0)--(1,e,0));
draw((1,0,-e)--(1,0,e));

triple c = (1,1,0);
path3 Q = plane_path3((1,-1,-1/2), c, theta=15, scale=2*e-0.3); 
drawing.Drawing3D.plane(Q,mediumred,nullpen,opacity=1);

path3 l = plane_intersection(Q,p); 
draw(l,black+1);

label("$\mathbb{T}^{2*}$", (1.5,-(e-0.5),e));
\end{asy}
\end{subfloatenv}
\begin{subfloatenv}{ }
\begin{asy}
import Drawing3D;
real e = 4;
DrawingR3 drawing = DrawingR3(e,0.1, axes_extent=e, camera=(10,2,7));
drawing.dual_axes();

path3 p = plane_path3((1,0,0),(1,0,0), v=(0,1,0), scale=2*e);
drawing.Drawing3D.plane(p);
draw((1,-e,0)--(1,e,0));
draw((1,0,-e)--(1,0,e));

triple c = (1,0,0);
path3 Q = plane_path3((0,-1,-1/2), c, theta=10, scale=2*e-0.6); 
drawing.Drawing3D.plane(Q,mediumred,nullpen,opacity=1);

path3 l = plane_intersection(Q,p); 
draw(l,black+1);

label("$\mathbb{T}^{2*}$", (1.5,-(e-0.5),e));
\end{asy}
\end{subfloatenv}\hfill%
\begin{subfloatenv}{ }
\begin{asy}
import Drawing3D;
real e = 4;
DrawingR3 drawing = DrawingR3(e,0.1, axes_extent=e, camera=(10,2,7));
drawing.dual_axes();

path3 p = plane_path3((1,0,0),(1,0,0), v=(0,1,0), scale=2*e);
drawing.Drawing3D.plane(p);
draw((1,-e,0)--(1,e,0));
draw((1,0,-e)--(1,0,e));

path3 l = line_path3((1,3,1)/4, direction=(1,3,1), extent_from_center=e);
draw(l,red+1);
dot((1,3,1));

draw((1,3,0)--(1,3,1), Dotted);
draw((1,0,1)--(1,3,1), Dotted);

label("$\mathbb{T}^{2*}$", (1.5,-(e-0.5),e));
\end{asy}
\end{subfloatenv}
\begin{subfloatenv}{ }
\begin{asy}
import Drawing3D;
real e = 4;
DrawingR3 drawing = DrawingR3(e,0.1, axes_extent=e, camera=(10,2,7));
drawing.dual_axes();

path3 p = plane_path3((1,0,0),(1,0,0), v=(0,1,0), scale=2*e);
drawing.Drawing3D.plane(p);
draw((1,-e,0)--(1,e,0));
draw((1,0,-e)--(1,0,e));

path3 l = line_path3((1,0,0)/4, direction=(1,0,0), extent_from_center=e);
draw(l,red+1.01);
dot((1,0,0));

label("$\mathbb{T}^{2*}$", (1.5,-(e-0.5),e));
\end{asy}
\end{subfloatenv}
\caption{Embedding of the dual space \T{2*} into \R{3*}}
\label{embedding T2* in R3*}
\end{figure}

The same constructions can be carried out for the dual space \T{2*}.
The space \R{3*} consists of the triples \((d,a,b)\),
where \((1,a,b)\in\R{3*}\) is identified with \((a,b)\in\T{2*}\),
so that a point \((a,b)\) can be identified with the one-dimensional linear subspace of \R{3*} that passes through \((1,a,b)\).
A plane in \R{3*} defined by \(d=1\) can be identified with the set of points of the dual space represented by sheaves.
The points represented by stacks are identified with one-dimensional linear subspaces which lie in the plane \(d=0\).
Lines are identified with two-dimensional linear subspaces.
In particular, the line that consists of the points represented by stacks is identified with the two-dimensional linear subspace \(d=0\).
Examples are shown in Figure~\ref{embedding T2* in R3*}, where the subspaces shown are identified with 
\(1-a-\tfrac{1}{2}b=0\) (a), \(2a+b=0\) (b), \((a,b)=(3,1)\) (c), and the origin (d).

Since points and lines are in one-to-one correspondence with linear subspaces,
the transformations \(\J\) and \(\Id\)
defined previously on points and lines of the dual space \T{2*} 
can be extended to one- and two-dimensional linear subspaces of \R{3*} in the obvious way.
In addition, I set \(\J(0,0,0)=\R{3}\) and \(\J(\R{3*})=(0,0,0)\) to complete the definition of \(\J\) and establish
a one-to-one correspondence between linear subspaces of \R{3} and \R{3*}.
I also set \(\Id(0,0,0)=(0,0,0)\) and \(\Id(\R{3*})=\R{3}\) to complete the definition of \(\Id\).
Observe that the subspaces shown in Figure~\ref{embedding T2* in R3*} are dual to the corresponding subspaces shown in Figure~\ref{embedding T2 in R3},
i.e.\ the one-dimensional subspace identified with \((-1,-\tfrac{1}{2})\in\T{2}\) as shown in Figure~\ref{embedding T2 in R3}(a)
and the two-dimensional subspace identified with the line \(1-a-\tfrac{1}{2}b=0\) as shown in Figure~\ref{embedding T2* in R3*}(a) are dual to each other,
and so on for the other panels.

Vectors of \R{3*} can be interpreted as linear functionals that act of vectors of the model space \R{3} and yield real numbers.
The value of the functional \((d,a,b)\in\R{3*}\) on \((w,x,y)\in\R{3}\) is defined by 
\begin{equation}
(d,a,b)[(w,x,y)]=dw+ax+by.
\end{equation}
The kernel of a functional \((d,a,b)\ne(0,0,0)\), which is defined as a set of all vectors \((w,x,y)\in\R{3}\) obeying \((d,a,b)[(w,x,y)]=0\),
is a two-dimensional linear subspace of \R{3} whose intersection with the plane \(w=1\) is a line given by \(d+ax+by=0\).
Substituting \(d=1\), I obtain the equation \(1+ax+by=0\) that determines projective duality between \T{2} and \T{2*}.
So, the space \R{3*} can be thought of as the usual vector-space dual of the model space \R{3}, where \R{3*} is defined as the space of the linear functionals acting on \R{3}.

This completes a brief introduction to projective geometry.
In the following sections, I focus on geometry of the target space \T{2}.
Every point or line in \T{2} can be represented by a specific linear subspace of \R{3}, with which it is identified as explained above.
Since  linear subspaces of \R{3*} and \R{3} are in one-to-one correspondence via the duality transformation,
every point or line in \T{2} can be represented by a specific linear subspace of \R{3*}, which is dual to the subspace of \R{3} identified with the point or the line.
The metric and the Clifford algebra are to be defined in \R{3*} rather than \R{3}, so
the latter indirect representation of points and lines in \T{2} via the duality transformation will be preferred.
Its further advantage is that it is consistent with the top-down model of \T{2}, which enables a more direct understanding of 
some geometric transformations as will become clear in the following.
The model space \R{3} is useful for certain non-metric aspects of geometry.
The projective dual \T{2*} is necessary for constructing \R{3*}, but it will only be used as an auxiliary space in the following.

\section{Grassmann algebra}
\label{section.grassmann}
\subsection{The basics}
Linear subspaces of \R{3} (and \R{3*}) and their intersections are described by a Grassmann algebra 
\cite{dorst2009geometric,roman2008advanced,greub1967multilinear}
also known as exterior algebra or a closely related Grassmann-Cayley algebra \cite{white1994grassmann}.
The Grassmann algebra of \R{3*} plays a more important role in the projective model since the metric 
(and Clifford algebra) are defined over \R{3*},
so I will consider it first.

The Grassmann algebra of \R{3*}  is denoted by \(\bigwedge\R{3*}\) and its elements are called multivectors.
It is an abstract eight-dimensional vector space (multivectors are ``vectors" of this abstract vector space),
whose basis consists of one scalar, three vectors, three bivectors, and one trivector, which is also called a pseudoscalar.
Scalars, vectors, bivectors, and trivectors are referred to as \(k\)-vectors where \(k=0\) for scalars, \(k=1\) for vectors, and so on
(the integer \(k\) is called the grade of a \(k\)-vector).
Addition of multivectors and multiplication by real numbers are defined in the usual componentwise way.
General multivectors are formed by linear combinations of scalars, vectors, bivectors, and trivectors.

Besides addition and multiplication by real numbers as in any vector space, 
it is also possible to compute the product of two multivectors, which turns \(\bigwedge\R{3*}\) into an algebra.
The product of multivectors in the Grassmann algebra is called the outer product (exterior product or wedge product are also used interchangeably).
The outer product of any two vectors \(\tb{a}, \tb{b}\in\R{3*}\) is anticommutative:
\begin{equation}
\tb{a}\wedge\tb{b}=-\tb{b}\wedge\tb{a},
\end{equation}
which implies \(\tb{a}\wedge\tb{a}=0\) for any \(\tb{a}\in\R{3*}\).
The outer product is neither commutative nor anticommutative in general,
but it is associative and distributive: 
\begin{equation}
\begin{split}
&A\wedge(B\wedge C)=(A\wedge B)\wedge C, \\
&A\wedge (B+ C)=A\wedge B+A\wedge C,\quad
(A+B)\wedge C=A\wedge C+B\wedge C
\end{split}
\end{equation}
for any \(A,B,C\in\bigwedge\R{3*}\).
Associativity allows one to drop the parentheses and write \(A\wedge B\wedge C\) without introducing an ambiguity.
Furthermore, \(s\wedge A=A\wedge s=s A\) for any scalar \(s\) and \(A\in\bigwedge\R{3*}\).
Given these properties, 
it is easy to verify that \(\tb{a}\wedge\tb{P}=\tb{P}\wedge\tb{a}\) for any vector \(\tb{a}\) and bivector \(\tb{P}\).

I let \(\e_0=(1,0,0)\), \(\e_1=(0,1,0)\), \(\e_2=(0,0,1)\) denote the standard basis vectors of \R{3*}, so that 
a vector \((d,a,b)\in\R{3*}\) can be written as \((d,a,b)=d\e_0+a\e_1+b\e_2\).
For brevity, I will use a simplified notation for the outer product of the basis vectors, e.g.\ \(\e_{01}=\e_0\wedge\e_1\) and \(\e_{012}=\e_0\wedge\e_1\wedge\e_2\).
The basis of the Grassmann algebra \(\bigwedge\R{3*}\) consists of the following \(2^3=8\) multivectors:
1, 
\(\e_0, \e_1, \e_2\), 
\(\e_{12}, \e_{20}, \e_{01}\),
\(\e_{012}\) (the basis trivector will also be denoted by \(\I=\e_{012}\)).
A general multivector in  \(\bigwedge\R{3*}\) can be written as
\begin{equation}
A=s+d\e_0+a\e_1+b\e_2+w\e_{12}+x\e_{20}+y\e_{01}+p\e_{012},
\end{equation}
where \(s,d,a,b,w,x,y,p\in\R{}\) can be thought of as the coordinates of \(A\)
in the abstract eight-dimensional vector space \(\bigwedge\R{3*}\).
The outer product of two multivectors can be computed by expressing both multivectors in terms of the basis 
of  \(\bigwedge\R{3*}\)
and
applying the properties of the outer product noted above, such as associativity, distributivity, 
and anticommutativity for the standard basis vectors \(\e_0, \e_1, \e_2\in\R{3*}\).

The outer product of a \(k\)-vector \(A_k\) and an \(l\)-vector \(B_l\) is a \((k+l)\)-vector, provided that \(k+l\le3\).
If \(k+l>3\), the outer product yields zero, i.e.\ a trivector is the highest grade \(k\)-vector that can be obtained in \(\bigwedge\R{3*}\).
A simple \(k\)-vector  (also called a blade or a decomposable \(k\)-vector) is a multivector which can be written as the outer product of \(k\) vectors
(scalars and vectors are simple by definition),
e.g.\ \(\e_{12}=\e_1\wedge\e_2\) is the outer product of two vectors and therefore it is a simple bivector.
In  \(\bigwedge\R{3*}\), all \(k\)-vectors including bivectors and trivectors 
are simple (an example of a non-simple bivector is \(\e_{01}+\e_{23}\) in \(\bigwedge\R{4*}\)).

The Grassmann algebra of the model space \R{3} has the same properties.
To avoid confusion, I will use \(\vee\) to denote the outer product in 
the Grassmann algebra of \R{3}, which will be denoted by \(\bigvee\R{3}\).
In speech it is convenient to refer to \(\vee\) as the join, because it corresponds to the join of subspaces as will 
be clarified in the following.
I let \(\e^0=(1,0,0)\), \(\e^1=(0,1,0)\), \(\e^2=(0,0,1)\) denote the standard basis vectors of \R{3}, 
so that a vector \((w,x,y)\in\R{3}\) can be written as \((w,x,y)=w\e^0+x\e^1+y\e^2\).
A simplified notation for bivectors and trivectors will be used in \(\bigvee\R{3}\), e.g.~\(\e^{12}=\e^1\vee\e^2\).
The basis of  \(\bigvee\R{3}\) consists of the following multivectors: 
1, 
\(\e^0, \e^1, \e^2\), 
\(\e^{12}, \e^{20}, \e^{01}\),
\(\e^{012}\).
The definitions given above for various concepts in \(\bigwedge\R{3*}\)  are adapted to \(\bigvee\R{3}\) in the obvious way.

\subsection{Blades and geometric objects}
A simple \(k\)-vector \(A_k\in\bigwedge\R{3*}\) represents a \(k\)-dimensional linear subspace of \R{3*}
consisting of all vectors \(\tb{a}\in\R{3*}\) which satisfy
\begin{equation}
\tb{a}\wedge A_k =0.
\label{subspace of blade}
\end{equation}
It can therefore represent a specific geometric object in \T{2*}, such as a point or a line,
and via projective duality can represent a specific geometric object in \T{2}.
The best way to understand this important point is to consider several examples given below.

The two-dimensional linear subspace shown in Figure~\ref{embedding T2* in R3*}(a) can be represented by the bivector \(\tb{P}=\e_{12}-\e_{20}-\tfrac{1}{2}\e_{01}\).
Indeed, substituting \(\tb{a}=d\e_0+a\e_1+b\e_2\) and \(A_k=\tb{P}\) into (\ref{subspace of blade}) gives
\begin{equation}
\notag
\tb{a}\wedge\tb{P}=(d\e_0+a\e_1+b\e_2)\wedge(\e_{12}-\e_{20}-\tfrac{1}{2}\e_{01})=(d-a-\tfrac{1}{2}b)\e_{012}=0,
\end{equation}
which is satisfied only if \(d-a-\tfrac{1}{2}b=0\).
The intersection of this two-dimensional linear subspace with the plane \(d=1\) is given by \(1-a-\tfrac{1}{2}b=0\),
which can be identified with a line in \T{2*}.
It is dual to the point \((-1,-\tfrac{1}{2})\) of the target space \T{2}.
Therefore, the bivector \(\tb{P}=\e_{12}-\e_{20}-\tfrac{1}{2}\e_{01}\) represents the point \((-1,-\tfrac{1}{2})\) shown in Figure~\ref{embedding T2 in R3}(a).
The bivector \(\tb{Q}=2\e_{20}+\e_{01}\) represents 
the two-dimensional linear subspace shown in Figure~\ref{embedding T2* in R3*}(b).
Substituting \(\tb{a}=d\e_0+a\e_1+b\e_2\) and \(A_k=\tb{Q}\) into (\ref{subspace of blade}) gives 
\begin{equation}
\notag
\tb{a}\wedge\tb{Q}=(d\e_0+a\e_1+b\e_2)\wedge(2\e_{20}+\e_{01}) =(2a+b)\e_{012}=0
\end{equation}
and thus \(2a+b=0\), which is a two-dimensional linear subspace of \R{3*} (as \(d\) ranges over \R{}).
Its intersection with the plane \(d=1\) is identified with the line in \T{2*} that passes through the origin of \T{2*}
and is dual to a point of the target space represented by the stack shown in Figure~\ref{embedding T2 in R3}(b).
The vector \(\tb{b}=\e_0+3\e_1+\e_2\) represents the one-dimensional linear subspace shown in Figure~\ref{embedding T2* in R3*}(c).
Substituting \(\tb{a}=d\e_0+a\e_1+b\e_2\) and \(A_k=\tb{b}\) into (\ref{subspace of blade}) gives 
\begin{equation}
\notag
\tb{a}\wedge\tb{b}=(d\e_0+a\e_1+b\e_2)\wedge(\e_0+3\e_1+\e_2)=
(a-3b)\e_{12}+(b-d)\e_{20}+(3d-a)\e_{01}=0
\end{equation}
and \(a=3b, b=d, a=3d\), where the first equation is redundant and the last two define a one-dimensional linear subspace,
which can be parametrised by \(d\) and consists of vectors \((d,3d,d)\).
A more general technique, which can be extended to higher-dimensional cases, is to write the three equations  in a matrix form as 
\begin{equation}
\notag
\begin{pmatrix}0&1&-3\\ 1&0&-1\\3&-1&0\end{pmatrix} \begin{pmatrix}d\\a\\b\end{pmatrix}=0.
\end{equation}
Since the rank of the matrix equals 2, e.g.\ the sum of the top and bottom rows can be expressed in terms of the middle row,
the kernel of the linear transformation defined by this matrix is a one-dimensional linear subspace of \R{3*}.
So, the equations \(a=3b, b=d, a=3d\) define a one-dimensional linear subspace and it is easy to verify that it passes through the point \((1,3,1)\in\R{3*}\)
by substituting \(d=1\). 
Hence, \(\tb{b}\) represents the point \((3,1)\) of the dual space \T{2*}
and via projective duality represents a line in the target space \T{2} defined by \(1+3x+y=0\).
Finally, substituting \(\tb{a}=d\e_0+a\e_1+b\e_2\) and \(A_k=\e_0\) into (\ref{subspace of blade}) gives 
\begin{equation}
\notag
\tb{a}\wedge\e_0=(d\e_0+a\e_1+b\e_2)\wedge \e_0=
b\e_{20}-a\e_{01}=0
\end{equation}
and thus \(a=0,b=0\), which is a one-dimensional linear subspace shown in Figure~\ref{embedding T2* in R3*}(d).
Its intersection with the plane \(d=1\) gives \((1,0,0)\) identified with the origin of \T{2*}, which is dual to the line in \T{2}
consisting of the points represented by stacks (see Figure~\ref{embedding T2 in R3}(d)).

In every example given above, I take the following sequence of steps:
\begin{equation}
\notag
\footnotesize
\text{a blade in }\bigwedge\R{3*}\to\text{a subspace of }\R{3*}\to\text{a geometric object of }\T{2*}\to\text{a geometric object of }\T{2},
\end{equation}
where the last step uses projective duality between \T{2} and \T{2*}.
It is also possible to arrive at the same geometric object in \T{2} by following the sequence:
\begin{equation}
\notag
\footnotesize
\text{a blade in }\bigwedge\R{3*}\to\text{a subspace of }\R{3*}\to\text{a subspace of }\R{3}\to\text{a geometric object of }\T{2},
\end{equation}
where the duality transformation defined on linear subspaces of \R{3*} is used in the second step.

In general, a bivector \(\tb{P}=\e_{12}+x\e_{20}+y\e_{01}\) represents a point in the target space located at \((x,y)\).
A bivector \(\tb{P}=\e_0\wedge(a\e_1+b\e_2)\) represents a point in \T{2}, 
which can only be visualised as a stack of lines,
e.g.\ this point is at infinity in Euclidean space.
It lies on the line \(ax+by=0\) and any line \(ax+by+d=0\), where \(d\in\R{}\).
A vector \(\tb{b}=d\e_0+a\e_1+b\e_2\) represents a line in the target space defined by \(ax+by+d=0\).
If \(d=0\), the line passes through the origin, and if \(a=b=0\) and \(d\ne0\), the line consists of the points represented by stacks,
e.g.\ it is the line at infinity in Euclidean space.
The correspondence between blades in \(\bigwedge\R{3*}\) and geometric objects in \T{2}
they represent is summarised in Table~\ref{from GR3* to T2}.

\begin{table}[t!]
\centering
\begin{tabular}{|c|c|}
\hline
blade in \(\bigwedge\R{3*}\) & geometric object in \T{2} \\
\hline
\(\e_{12}+x\e_{20}+y\e_{01}\) & the sheaf of lines passing through \((x,y)\) \\
\(\e_0\wedge(a\e_1+b\e_2)\) & the stack of lines \(ax+by+d=0\), \(d\in\R{}\) \\
\(d\e_0+a\e_1+b\e_2\)  & the line \(ax+by+d=0\) if \((a,b)\ne(0,0)\) \\
\(\e_0\) & the line that consists of points represented by stacks \\
\(\e_1\) & the \(y\)-axis (\(x=0\))\\
\(\e_2\) & the \(x\)-axis (\(y=0\))\\
\hline
\end{tabular}
\caption{The correspondence between blades in \(\bigwedge\R{3*}\) and geometric objects in \T{2}}
\label{from GR3* to T2}
\end{table}

Since the subspace represented by a blade \(s A_k\), where the scalar \(s\ne0\), is exactly the same as the subspace represented by \(A_k\),
both \(s A_k\) and \(A_k\) represent the same geometric object in the target space.
A blade is not uniquely determined by its subspace; one needs to specify orientation and weight to fully determine a blade.
Both orientation and weight are relative concepts determined in relation to some other blade whose subspace is the same.
For instance, the vector \(\tb{a}=-2(\e_{0}+2\e_1)\) in relation to \(\tb{b}=\e_{0}+2\e_1\) has the opposite orientation and the weight of \(2\).
It is convenient to transfer the orientation and weight of a blade \(s A_k\) to the geometric object in the target space it represents.
I can then say that the vectors \(\tb{a}\) and \(\tb{b}\) represent the same line but with the opposite orientations and
that the weight of the line represented by \(\tb{a}\) is twice that of the line represented by \(\tb{b}\).
In this sense, blades can be identified with the geometric objects of the target space they represent.
Therefore, I can refer to a blade as a point or a line, e.g.\ \(\e_{12}\) is the origin of \T{2} and
\(\e_{2}\) is the \(x\)-axis.

The correspondence between linear subspaces of the model space \R{3}
and the \(k\)-vectors of the Grassmann algebra \(\bigvee\R{3}\) is established in the same way as in \(\bigwedge\R{3*}\).
Namely, a simple \(k\)-vector \(X_k\in\bigvee\R{3}\) represents a \(k\)-dimensional linear subspace of \R{3}
consisting of all vectors \(\tb{x}\in\R{3}\) which satisfy
\begin{equation}
\tb{x}\vee X_k =0.
\label{subspace of blade in R3}
\end{equation}
Therefore, a blade \(X_k\in\bigvee\R{3}\) represents a geometric object in the target space.
This representation is direct in contrast to the indirect representation considered previously, which relies on duality.
For example, a one-dimensional linear subspace shown in Figure~\ref{embedding T2 in R3}(a) 
can be represented by the vector 
\(\tb{u}=\e^0-\e^1-\tfrac{1}{2}\e^2\).
Substituting \(\tb{x}=w\e^0+x\e^1+y\e^2\) and \(X_k=\tb{u}\) into (\ref{subspace of blade in R3})
gives
\begin{equation}
\notag
\tb{x}\vee\tb{u}=(w\e^0+x\e^1+y\e^2)\vee(\e^0-\e^1-\tfrac{1}{2}\e^2)=
(y-\tfrac{1}{2}x)\e^{12}+(y+\tfrac{1}{2}w)\e^{20}-(x+w)\e^{01}=0
\end{equation}
and therefore \(y=\tfrac{1}{2}x\), \(y=-\tfrac{1}{2}w\), \(x=-w\), where
the first equation is redundant and the last two define the subspace of \R{3} that
consists of vectors \((w,-w,-\tfrac{1}{2}w)\).
This subspace intersects the plane \(w=1\) at the point \((1,-1,-\tfrac{1}{2})\) identified with 
the point \((-1,-\tfrac{1}{2})\) of the target space.
So, \(\tb{u}=\e^0-\e^1-\tfrac{1}{2}\e^2\) represents the point \((-1,-\tfrac{1}{2})\).
Substituting \(\tb{x}=w\e^0+x\e^1+y\e^2\) and \(\tb{F}=\e^{12}+3\e^{20}+\e^{01}\) 
into (\ref{subspace of blade in R3}) gives
\begin{equation}
\notag
\tb{x}\vee \tb{F}=(w\e^0+x\e^1+y\e^2)\vee (\e^{12}+3\e^{20}+\e^{01})=(w+3x+y)\e^{012}=0
\end{equation}
and thus \(w+3x+y=0\), which is a two-dimensional linear subspace shown in Figure~\ref{embedding T2 in R3}(c). 
It intersects the plane \(w=1\) along a line, which can be identified with a line in \T{2} defined by \(1+3x+y=0\).
So, the bivector \(\tb{F}\) represents a line  in \T{2} defined by \(1+3x+y=0\).

\begin{table}[t!]
\centering
\begin{tabular}{|c|c|}
\hline
blade in \(\bigvee\R{3}\) & geometric object in \T{2} \\
\hline \\[-10pt]
\(\e^0+x\e^1+y\e^2\)  & the point at  \((x,y)\) \\
\(d\e^{12}+a\e^{20}+b\e^{01}\) & the set of points on the line \(ax+by+d=0\) if \((a,b)\ne(0,0)\) \\
\hline
\end{tabular}
\caption{The correspondence between blades in \(\bigvee\R{3}\) and geometric objects in \T{2}}
\label{from GR3 to T2}
\end{table}

In general, a vector \(\tb{u}=\e^0+x\e^1+y\e^2\) represents a point in \T{2} located at \((x,y)\)
and a bivector \(\tb{F}=d\e^{12}+a\e^{20}+b\e^{01}\) represents a line in \T{2}, which
consists of the points satisfying \(ax+by+d=0\).
Vectors of the model space \R{3} confined to the plane \(w=0\) represent points, 
which in the top-down model are described by stacks of lines.
There is no obvious metric-independent way to refer to them in the bottom-up model of \T{2},
but if the space is Euclidean they can be described as points at infinity.
The line represented by \(\e^{12}\) is at infinity in Euclidean space as well.
The correspondence between blades in \(\bigvee\R{3}\) and geometric objects in \T{2}
is summarised in Table~\ref{from GR3 to T2}.
The blades \(X_k\) and \(s X_k\), where \(s\ne0\), represent the same linear subspace of \(\R{3}\)
and therefore they also represent the same geometric object in \T{2}.
It is convenient to transfer the orientation and weight of a blade \(s X_k\) to the geometric object it represents.
For example,  \(\tb{v}=-2(\e^0-\e^1-\tfrac{1}{2}\e^2)\) and \(\tb{u}=\e^0-\e^1-\tfrac{1}{2}\e^2\)
represent the same point located at \((-1,-\tfrac{1}{2})\) but with the opposite orientations,
and the weight of the point represented by \(\tb{v}\) is twice that of the point represented by \(\tb{u}\).
I will identify blades in \(\bigvee\R{3}\) with the geometric objects they represent in \T{2},
so that I can refer to a vector \(\tb{u}\) as a point and a bivector \(\tb{F}\) as a line.
To distinguish the orientation and weight due to the blades of \(\bigvee\R{3}\)
from those due to the blades of \(\bigwedge\R{3*}\), I will refer to the former
as the bottom-up orientation and weight and to the latter as the top-down orientation and weight.

\subsection{Meet and join}
The Grassmann algebra \(\bigwedge\R{3*}\) is related to the top-down model of \T{2}
and is therefore concerned with the intersection of lines as a way to represent points.
It also \ looks at the incidence of points and lines from the top-down point of view.
The subspace of \R{3*} represented by a vector \(\tb{a}\) lies in the subspace 
represented by a bivector \(\tb{P}\) if and only if \(\tb{a}\wedge\tb{P}=0\),
which implies that the line in the target space represented by \(\tb{a}\) belongs to 
the pencil of lines represented by \(\tb{P}\).
By identifying blades with the geometric objects as explained previously,
it is meaningful to say that a line \(\tb{a}\) passes through a point \(\tb{P}\) if and only if \(\tb{a}\wedge\tb{P}=0\).
The outer product \(\tb{a}\wedge\tb{b}\) of two lines \(\tb{a}\) and \(\tb{b}\) 
is a point where the lines intersect.
It corresponds to the meet of the two-dimensional subspaces of \R{3} identified with the lines.
As expected, both lines \(\tb{a}\) and \(\tb{b}\) pass through the point \(\tb{a}\wedge\tb{b}\),
since \(\tb{a}\wedge(\tb{a}\wedge\tb{b})=0\) and \(\tb{b}\wedge(\tb{a}\wedge\tb{b})=0\).

On the other hand, the Grassmann algebra \(\bigvee\R{3}\) is related to the bottom-up model of \T{2}.
It is therefore concerned with the join of points as a way to build lines
and it looks at the incidence of points and lines from the bottom-up point of view.
The subspace of \R{3} represented by a vector \(\tb{u}\in\R{3}\) lies in the subspace represented by 
a bivector \(\tb{F}\in\bigvee\R{3}\) if and only if \(\tb{u}\vee\tb{F}=0\),
which implies that the point \(\tb{u}\) lies on the line  \(\tb{F}\).
The outer product \(\tb{u}\vee\tb{v}\) of two points \(\tb{u}\) and \(\tb{v}\)  
is a line passing through the points.
It corresponds to the join of the one-dimensional subspaces of \R{3} identified with the points.
As expected, both points \(\tb{u}\) and \(\tb{v}\)  lie on the line \(\tb{u}\vee\tb{v}\), since
\(\tb{u}\vee(\tb{u}\vee\tb{v})=0\) and \(\tb{v}\vee(\tb{u}\vee\tb{v})=0\).

\begin{table}
\centering
\begin{tabular}{|l|cccccccc|}
\hline
\(A\) &         \(1\)       & \(\e_0\)&    \(\e_1\)&         \(\e_2\)&    \(\e_{12}\)& \(\e_{20}\)& \(\e_{01}\)& \(\e_{012}\)\\
\hline \\ [-10pt]
\(\J(A)\phantom{^{-1}}\) & \(-\e^{012}\) & \(-\e^{12}\)& \(-\e^{20}\)&  \(-\e^{01}\)&    \(\phantom{-}\e^0\)&    \(\phantom{-}\e^{1}\)& \(\phantom{-}\e^{2}\)& \(\phantom{-}1\) \\
\hline
\end{tabular}
\caption{The list of values of  \(\J:\bigwedge\R{3*}\to\bigvee\R{3}\) on the basis  multivectors of \(\bigwedge\R{3*}\).}
\label{duality R3}
\end{table}

It follows from the above discussion that one uses \(\bigwedge\R{3*}\) to compute the meet
and \(\bigvee\R{3}\) to compute the join.
In order to define the join in \(\bigwedge\R{3*}\) and the meet in \(\bigvee\R{3}\), 
one needs to extend the duality transformation to multivectors.
The definition of \(\J:\bigwedge\R{3*}\to\bigvee\R{3}\) is given in \ref{appendix.duality Rn}. 
The values \(\J\) assumes on the basis multivectors of \(\bigwedge\R{3*}\) are listed in Table~\ref{duality R3}.
This table and the property of linearity can be used to compute \(\J(A)\) for any multivector \(A\in\bigwedge\R{3*}\), e.g.\ 
\(\J(\e_{12}+x\e_{20}+\e_{01})=\e^0+x\e^1+y\e^2\) and \(\J(d\e_0+a\e_1+b\e_2)=-(d\e^{12}+a\e^{20}+b\e^{01})\).
The inverse \(\J^{-1}:\bigvee\R{3}\to\bigwedge\R{3*}\) can be computed using the same table.

For \(A, B\in\bigwedge\R{3*}\),  the join of \(A\) and \(B\) is a multivector in \(\bigwedge\R{3*}\) defined by
\begin{equation}
A\vee B = \J^{-1}(\J(A)\vee\J(B)),
\end{equation}
where I use the same symbol \(\vee\) for the join in \(\bigwedge\R{3*}\) as the one used for the outer product in \(\bigvee\R{3}\).
The join is metric-independent, since it only depends on the duality transformation and the outer product in \(\bigvee\R{3}\).
The join thus defined allows one to gain access to the bottom-up model of \T{2} without leaving \(\bigwedge\R{3*}\).
For \(F,G\in\bigvee\R{3}\), the meet of \(F\) and \(G\) is a multivector in \(\bigvee\R{3}\) defined by 
\(F\wedge G = \J(\J^{-1}(F)\wedge\J^{-1}(G))\),
where I use the same symbol \(\wedge\) for the meet in \(\bigvee\R{3}\) as the one used for the outer product in \(\bigwedge\R{3*}\).
I will not find much need for the meet in \(\bigvee\R{3}\), since 
any non-metric problem can be addressed in  \(\bigwedge\R{3*}\) by means of \(\wedge\) and \(\vee\) defined there
and \(\bigwedge\R{3*}\) will be the focus of the ensuing discussion when the metric is introduced.

The join \(\tb{P}\vee\tb{Q}\) of two points \(\tb{P}\) and \(\tb{Q}\) is
a line passing through the points. It satisfies \(\tb{Q}\vee\tb{P}=-\tb{P}\vee\tb{Q}\),
i.e.\ the join of two points is anticommutative, which implies that
the bottom-up orientation of the line \(\tb{Q}\vee\tb{P}\) is the opposite of that of the line \(\tb{P}\vee\tb{Q}\).
The equality 
\begin{equation}
\tb{a}\wedge\tb{P}=(\tb{a}\vee\tb{P})\I,
\end{equation}
is satisfied for any line \(\tb{a}\) and point \(\tb{P}\).
One can easily verify that  the join of a line and a point are commutative, i.e.\ 
 \(\tb{a}\vee\tb{P}=\tb{P}\vee\tb{a}\).
The join of a multivector \(A\) with the basis pseudoscalar \(\I\) yields the same multivector, \(A\vee\I=\I\vee A=A\).
Furthermore, \(A_k\vee B_l=0\) if \(k+l<3\), where \(A_k\) is a \(k\)-vector and \(B_l\) is an \(l\)-vector, e.g.\ 
\(\tb{a}\vee\tb{b}=0\) and \(s\vee\tb{P}=0\).

I illustrate the use of the join by solving the problem of finding the central point of a line defined by \(ax+by+1=0\).
To find the central point, I first  find 
the line that passes through the origin \(\tb{P}=\e_{12}\) and the point \(\tb{Q}=\e_{12}+a\e_{20}+b\e_{01}\).
It is given by
\begin{multline*}
\tb{P}\vee\tb{Q}=\e_{12}\vee(\e_{12}+a\e_{20}+b\e_{01})=
\J^{-1}(\J(\e_{12})\vee\J(\e_{12}+a\e_{20}+b\e_{01}))=\\
=\J^{-1}(\e^0\vee(\e^0+a\e^1+b\e^2))=\J^{-1}(-b\e^{20}+a\e^{01})=b\e_1-a\e_2.
\end{multline*}
Then the intersection of the lines \(\tb{a}=\e_0+a\e_1+b\e_2\) and \(\tb{P}\vee\tb{Q}=b\e_1-a\e_2\) is given by
\begin{equation}
\notag
\tb{a}\wedge(\tb{P}\vee\tb{Q})=(\e_0+a\e_1+b\e_2)\wedge(b\e_1-a\e_2)=
-(a^2+b^2)\e_{12}+a\e_{20}+b\e_{01},
\end{equation}
which yields the coordinates of the central point of \(\tb{a}\).

\section{The metric and Clifford algebra}
\label{section.clifford}
\subsection{The basics}
The Clifford algebra of \R{3*} deals with the orthogonality relations between the subspaces of \R{3*}
and it cannot be defined without specifying a metric.
I assume the standard basis vectors  are orthogonal, i.e.\ \(\e_i\cdot\e_j=0\) for \(i\ne j\).
The inner product of vectors \(\tb{a}\) and \(\tb{b}\)  is given by
\begin{equation}
\tb{a}\cdot\tb{b}=d_1d_2\e_0\cdot\e_0+a_1a_2\e_1\cdot\e_1+b_1b_2\e_2\cdot\e_2,
\end{equation}
where \(\tb{a}=d_1\e_0+a_1\e_1+b_1\e_2\) and  \(\tb{b}=d_2\e_0+a_2\e_1+b_2\e_2\).
The metric is specified by assigning values to the inner product of the standard basis vectors with themselves: 
\(\e_0\cdot\e_0, \e_1\cdot\e_1, \e_2\cdot\e_2\).
Some of the choices of interest are given in Table~\ref{metric}.
Each choice of the metric engenders specific metric properties in the target space \T{2},
which becomes a homogeneous metric space;
the metric spaces defined in the projective model are known as Cayley-Klein geometries.

\begin{table}
\centering
\begin{tabular}{|ccc|c|l|c|c|}
\hline
\(\e_0\cdot\e_0\) & \(\e_1\cdot\e_1\) & \(\e_2\cdot\e_2\)  & \T{2} & space & distances & angles\\
\hline
0 & 1 & 1 &\E{2}& Euclidean & degenerate & elliptic\\ 
 1& 1& 1& \El{2} & Elliptic & elliptic & elliptic\\
\m1 & 1& 1 & \Hy{2} & Hyperbolic & hyperbolic & elliptic\\
0 & 1  & \m1 & \M{2} &Minkowski  & degenerate & hyperbolic\\
 \m1 & 1 & \m1 & \AdS{2} & Anti de-Sitter & elliptic & hyperbolic\\
1 & 1 & \m1 & \dS{2} & de-Sitter & hyperbolic & hyperbolic \\
\hline
\end{tabular}
\caption{Some choices of the metric and the corresponding homogeneous metric spaces.}
\label{metric}
\end{table}

The metric is called degenerate if \(\e_i\cdot\e_i=0\) for some \(i\); the corresponding metric space is also called degenerate.
The sign of  \(\e_0\cdot\e_0\) is indicative of the curvature of the corresponding space and 
\(\e_0\cdot\e_0=0\) implies that the space is flat (zero curvature).
Euclidean and Minkowski spaces are flat and degenerate; the other spaces listed in Table~\ref{metric} are non-degenerate.
Elliptic space is closely related to spherical geometry 
(without going into details, I note that elliptic space is equivalent to spherical space, where the antipodal points are identified).
Minkowski, de-Sitter, and anti de-Sitter spaces admit a relativistic kinematic structure.
I will refer to Minkowski, de-Sitter, and anti de-Sitter spaces as kinematic and 
Euclidean, elliptic, and hyperbolic spaces as non-kinematic.
There are nine possible Cayley-Klein geometries in two dimensions.
Besides the six listed in Table~\ref{metric}, there are three geometries with degenerate angle measure.
These geometries will not be considered in this paper.
The following discussion is limited to the six geometries listed in Table~\ref{metric}.

The Clifford algebra of \R{3*} is closely related to the Grassmann algebra of \R{3*}.
In fact, as abstract vector spaces they are identical.
Namely, the Clifford algebra of \R{3*} consists of the same multivectors as \(\bigwedge\R{3*}\),
its basis multivectors are the same as  in  \(\bigwedge\R{3*}\),
and multivector addition and multiplication by real numbers are the same as in \(\bigwedge\R{3*}\).
What makes it distinct from the Grassmann algebra is the product of multivectors,
which is called the geometric (or Clifford) product. It is denoted by juxtaposition of multivectors
(the geometric product of a multivector \(A\) with itself is usually denoted by \(A^2=AA\)).
The geometric product of any two vectors \(\tb{a}, \tb{b}\in\R{3*}\) yields a multivector with a scalar and a bivector component:
\begin{equation}
\tb{a}\tb{b}=\tb{a}\cdot\tb{b}+\tb{a}\wedge\tb{b},
\end{equation}
which implies \(\e_i^2=\e_i\cdot\e_i\), \(i=0,1,2\), and \(\e_i\e_j=-\e_j\e_i\) if \(i\ne j\).
The geometric product is associative and distributive:
\begin{equation}
\begin{split}
&A(B C)=(AB) C, \\
&A (B+ C)=A B+A C,\quad
(A+B) C=A C+B C
\end{split}
\end{equation}
for any \(A,B,C\in\bigwedge\R{3*}\).
Since the geometric product is associative, one can drop the parentheses and write \(ABC\).
Furthermore, the geometric product \(sA=As\) of a scalar \(s\) and a multivector \(A\) is the same as multiplication by \(s\)
in the abstract vector space, i.e.\ it is defined componentwise 
(since both operations are denoted by juxtaposition, it is not even possible to distinguish between them in the standard notation).
Since the basis vectors of \R{3*} are assumed to be orthogonal, the geometric and outer products are identical on the basis multivectors,
e.g.\ \(\e_0\e_1\e_2=\e_0\wedge\e_1\wedge\e_2\).
The geometric product of any multivectors \(A\) and \(B\) can be computed by expressing both multivectors in terms of the basis multivectors
and applying the properties noted above.

When using a Clifford algebra, it is not uncommon to deal with multivectors consisting of several components with different grades.
I let \(\grade{A}_k\) denote the \(k\)-vector component of a multivector \(A\).
The geometric product of a \(k\)-vector \(A_k\) and an \(l\)-vector \(B_l\) yields a multivector 
with two components, which have the grades \(|k-l|\) and \(|k-l|+2\) (the latter only if it is \(\le3\)).
The outer product of \(A_k\) and \(B_l\) is inherited from the Grassmann algebra and satisfies 
\begin{equation}
A_k\wedge B_l=\grade{A_kB_l}_{k+l}.
\end{equation}
The inner product (also called the dot product) of \(A_k\) and \(B_l\) is defined by
\begin{equation}
A_k\cdot B_l=\grade{A_kB_l}_{|k-l|},
\end{equation}
which extends by linearity to general multivectors.
The basis pseudoscalar \(\I\) of \(\bigwedge\R{3*}\) commutes with any multivector \(A\), i.e.\ \(\I A=A\I\),
and \(A\I=A\cdot\I\).
The inner product of a \(k\)-vector \(A_k\) and an \(l\)-vector \(B_l\), where \(k\le l\),  is related to the join and the outer product by
\begin{equation}
A_k\cdot B_l = (A_k\I)\vee B_l \quad \text{and} \quad
(A_k\cdot B_l)\I = A_k\wedge (B_l\I).
\label{dot to vee}
\end{equation}
In particular, \(\tb{a}\cdot\tb{b}=(\tb{a}\I)\vee\tb{b}\) and \((\tb{a}\cdot\tb{P})\I=\tb{a}\wedge(\tb{P}\I)\),
where \(\tb{a}\) is a line and \(\tb{P}\) is a point. 
The commutator of two multivectors \(A\) and \(B\) is defined by
\begin{equation}
A\times B=\tfrac{1}{2}(AB-BA),
\end{equation}
which implies \(A\times B=-B\times A\) for any \(A,B\in\bigwedge\R{3*}\).
The commutator satisfies the following identities: 
\begin{equation}
\tb{a}\times\tb{b}=\tb{a}\wedge\tb{b},\quad
\tb{a}\times\tb{P}=\tb{a}\cdot\tb{P},\quad \text{and}\quad 
\tb{P}\times\tb{Q}=(\tb{P}\vee\tb{Q})\I,
\end{equation}
where \(\tb{a}\), \(\tb{b}\) are lines and \(\tb{P}\), \(\tb{Q}\) are points.
It follows that the commutator of two points is again a point.
The properties of the geometric product can be summarised as follows:
\begin{equation}
\begin{split}
&\tb{a}\tb{b}=\tb{a}\cdot\tb{b}+\tb{a}\wedge\tb{b},\quad \text{where}\quad
\tb{a}\cdot\tb{b}=\tb{b}\cdot\tb{a},\quad
\tb{a}\wedge\tb{b}=-\tb{b}\wedge\tb{a},\\
&\tb{a}\tb{P}=\tb{a}\cdot\tb{P}+\tb{a}\wedge\tb{P},\quad \text{where} \quad
\tb{a}\cdot\tb{P}=-\tb{P}\cdot\tb{a},\quad
\tb{a}\wedge\tb{P}=\tb{P}\wedge\tb{a},\\
&\tb{P}\tb{Q}=\tb{P}\cdot\tb{Q}+(\tb{P}\vee\tb{Q})\I,\quad \text{where} \quad
\tb{P}\cdot\tb{Q}=\tb{Q}\cdot\tb{P}, \quad
\tb{P}\vee\tb{Q}=-\tb{Q}\vee\tb{P}.
\end{split}
\end{equation}
The reverse \(\reverse{A}\) of a multivector \(A=s+\tb{a}+\tb{P}+p\I\),
where \(\tb{a}\) is a line and \(\tb{P}\) is a point, is given by
\begin{equation}
\reverse{A}=s+\tb{a}-\tb{P}-p\I.
\end{equation}
It satisfies \(\reverse{AB}=\reverse{B}\reverse{A}\)  for any \(A,B\in\bigwedge\R{3*}\).
I let \((A)_G\) denote a multivector obtained from \(A\) by reversing the sign of all \(k\)-vector components of \(A\) for which the grade \(k\)
is in the list of grades \(G\), e.g.\ \(\reverse{A}=(A)_{23}\).
The norm of a multivector \(A\) is defined by
\begin{equation}
\norm{A}=|A\reverse{A}(A\reverse{A})_1|^{\frac{1}{4}},
\end{equation}
where \(A\reverse{A}(A\reverse{A})_1\) is a scalar for any  \(A\in\bigwedge\R{3*}\) \cite{dadbeh2011inverse,shirokov2011concepts}.
It satisfies \(\norm{AB}=\norm{A}\norm{B}\) for any \(A,B\in\bigwedge\R{3*}\), e.g.\ \(\norm{sA}=|s|\norm{A}\) where \(s\) is a scalar.
If \(A\) is a blade, then \(A\reverse{A}\) is a scalar and \(\norm{A}=|A\reverse{A}|^{\frac{1}{2}}\).
A multivector \(A\) is called normalised if \(\norm{A}=1\). Any multivector whose norm is not zero can be normalised.
Note that I use the term \emph{norm} informally, since \(\norm{A}=0\) does not necessarily imply \(A=0\) and 
\(\norm{A+B}\le\norm{A}+\norm{B}\) may not be satisfied in some spaces.
If the norm of \(A\) is not zero, \(A\) has the inverse given by
\begin{equation}
A^{-1}=
\frac{\reverse{A}(A\reverse{A})_1}{A\reverse{A}(A\reverse{A})_1}.
\end{equation}
If \(A\) is a blade with the nonzero norm, then the inverse is given by \(A^{-1}=\reverse{A} / (A\reverse{A})\).
The exponential function is defined on an arbitrary multivector \(A\) by the Taylor series
\begin{equation}
e^A=\sum_{n=0}^{\infty}\frac{A^n}{n!}.
\end{equation}
For commuting multivectors \(A,B\in\bigwedge\R{3*}\), i.e.~\(AB=BA\), the equality \(e^{A+B}=e^Ae^B\) holds.
The inverse of \(e^A\) is given by \(e^{-A}\) for any \(A\in\bigwedge\R{3*}\).
If \(A\) is a blade with the zero norm, then \(A^2=0\) and therefore \(e^A=1+A\).
If \(A\) is a blade with the nonzero norm, it can be normalised and therefore it is sufficient to consider \(e^{\alpha A}\), where
\(\alpha\) is a scalar and \(A^2=1\) or \(A^2=-1\), which gives
\begin{equation}
e^{\alpha A}=\left\{
\begin{split}
&1+\alpha A, \text{ if }A^2=0,\\
&\cos\alpha+A\sin\alpha, \text{ if }A^2=-1,\\
&\cosh\alpha+A\sinh\alpha, \text{ if }A^2=1.
\end{split}
\right.
\end{equation}
A multivector \(A\) is called even if \(\grade{A}_k=0\) for all odd grades \(k\).
An even multivector can be written as
\begin{equation}
A=s+w\e_{12}+x\e_{20}+y\e_{01}.
\end{equation}
Even multivectors form a subalgebra of the Clifford algebra, since the geometric product of two even multivectors is even.
If the metric is elliptic, the subalgebra of even multivectors is isomorphic to the algebra of quaternions.

\subsection{Distances and angles}
Selecting a metric in \R{3*} makes it possible to measure distances and angles in the target space \T{2}.
The procedure for determining the distance between a pair of points or the angle between a pair of lines
employed in Cayley-Klein geometries is quite involved as it relies on synthetic reasoning.
The Clifford algebra of \R{3*} simplifies this procedure by replacing geometric reasoning with algebraic formulas.

The projective model provides access to three distinct kinds of geometric objects: proper,  improper, and null.
The metric content of the target space \T{2} is restricted to proper geometric objects.
In non-kinematic spaces, a blade \(A\) is called proper if \(A\reverse{A}>0\), improper if \(A\reverse{A}<0\), and null otherwise.
In kinematic spaces, the same applies to all blades except points, for which the opposite sign is required,
i.e.\ a point \(\tb{P}\) in a kinematic space is called proper if \(\tb{P}\reverse{\tb{P}}<0\) and improper if \(\tb{P}\reverse{\tb{P}}>0\)
(the change of sign is necessary, since I want the origin to be a proper point 
but \(\e_{12}\reverse{\e}_{12}=\e_1^2\e_2^2=-1\) in kinematic spaces).
The distance \(r\) between two \emph{proper} points \(\tb{P}\) and \(\tb{Q}\) 
is defined if the join \(\tb{P}\vee\tb{Q}\) is a \emph{proper} line.
If \(\tb{P}\) and \(\tb{Q}\) are \emph{normalised}, the distance is determined by
\begin{equation}
m_d(r)=\norm{\tb{P}\vee\tb{Q}}
\end{equation}
where \(m_d(r)\) depends on the distance measure of the metric space (see Table~\ref{metric}) and is given by
\begin{equation}
m_d(r)=\left\{
\begin{split}
&r, \text{ if distance measure is degenerate (Euclidean, Minkowski)},\\
&\sin r, \text{ if distance measure is elliptic (Elliptic, Anti de-Sitter)},\\
&\sinh r, \text{ if distance measure is hyperbolic (Hyperbolic, de-Sitter)}.
\end{split}\right.
\end{equation}
The distance \(r\) ranges from 0 to \(+\infty\) if the distance measure is degenerate or hyperbolic,
and it ranges from 0 to \(\tfrac{\pi}{2}\) if the distance measure is elliptic.
Hence,  the maximum distance between points equals \(\tfrac{\pi}{2}\) in elliptic and anti de-Sitter spaces.
Note that in non-kinematic spaces, the line  \(\tb{P}\vee\tb{Q}\) is proper if the points \(\tb{P}\) and \(\tb{Q}\) are proper,
so the distance is defined between any proper points.
The angle \(\alpha\) between two \emph{proper} lines \(\tb{a}\) and \(\tb{b}\) is defined 
if their intersection  \(\tb{a}\wedge\tb{b}\) is a \emph{proper} point.
If \(\tb{a}\) and \(\tb{b}\) are \emph{normalised}, the angle is determined by
\begin{equation}
m_a(\alpha)=|\tb{a}\cdot\tb{b}|
\end{equation}
where \(m_a(\alpha)\) depends on the angle measure of the metric space and is given by
\begin{equation}
m_a(\alpha)=\left\{
\begin{split}
&\cos\alpha, \text{ if angle measure is elliptic (non-kinematic spaces)},\\
&\cosh\alpha, \text{ if angle measure is hyperbolic (kinematic spaces)}.
\end{split}\right.
\end{equation}
The angle \(\alpha\) ranges from 0 to \(\tfrac{\pi}{2}\) in non-kinematic spaces, and from 0 to \(+\infty\) in kinematic spaces.
In non-kinematic spaces, one can use \(\cos\alpha=\tb{a}\cdot\tb{b}\) where \(\alpha\in[0,\pi]\)
if it is desirable to take the orientation of the lines into account.
The lines \(\tb{a}\) and \(\tb{b}\) are called perpendicular if \(\tb{a}\cdot\tb{b}=0\), 
which gives \(\alpha=\tfrac{\pi}{2}\) in non-kinematic spaces (the angle between perpendicular lines is undefined in kinematic spaces).
Since \(\tb{a}\cdot\tb{b}=(\tb{a}\I)\vee\tb{b}\) for any lines \(\tb{a}\) and \(\tb{b}\), 
a line \(\tb{b}\) is perpendicular to \(\tb{a}\) if and only if
it passes through the point \(\tb{a}\I\), which is called the polar point of the line \(\tb{a}\).
Note that a null line is perpendicular to itself according to the definition.
In non-degenerate spaces, if a line \(\tb{a}\) is null, then \(\tb{a}\I\) is a null point on \(\tb{a}\),
and if a point \(\tb{P}\) is null, then \(\tb{P}\I\) is a null line through \(\tb{P}\).

The inner product \(\tb{a}\cdot\tb{P}\) of a line \(\tb{a}\) and a point \(\tb{P}\) yields a line, 
which passes through  \(\tb{P}\) and is perpendicular to \(\tb{a}\), i.e.\ \(\tb{P}\wedge(\tb{a}\cdot\tb{P})=0\) and  \(\tb{a}\cdot(\tb{a}\cdot\tb{P})=0\).
For instance, to find a line in \M{2}, which is perpendicular to the line \(\tb{a}=\e_0+3\e_1+\e_2\) and passes through the origin \(\tb{P}=\e_{12}\), 
I compute \(\tb{a}\cdot\tb{P}=\grade{\tb{a}\tb{P}}_{|1-2|}=\grade{(\e_0+3\e_1+\e_2)\e_{12}}_1=
\grade{\e_{012}+3\e_1^2\e_2-\e_2^2\e_1}_1=\e_1+3\e_2\), which gives the line \(x+3y=0\).
The distance \(r\) between \tb{a} and \(\tb{P}\) is determined by \(m_d(r)=|\tb{a}\vee\tb{P}|\) 
if both \tb{a} and \(\tb{P}\) are \emph{normalised}
and if the line \(\tb{a}\cdot\tb{P}\) and the points \(\tb{P}\) and \(\tb{a}\wedge(\tb{a}\cdot\tb{P})\) are \emph{proper}
(the distance from \(\tb{P}\) to \(\tb{a}\wedge(\tb{a}\cdot\tb{P})\) also equals \(r\)).
The inner product \(\tb{a}\cdot\tb{P}\) can be used in non-degenerate spaces to compute \(r\) thanks to the relations
\(\norm{\tb{a}\cdot\tb{P}}^2+|\tb{a}\vee\tb{P}|^2=1\) if the distance measure is elliptic
and \(\norm{\tb{a}\cdot\tb{P}}^2-|\tb{a}\vee\tb{P}|^2=1\) if the distance measure is hyperbolic.
Similarly, the inner product \(\tb{P}\cdot\tb{Q}\) can be used in non-degenerate spaces 
to compute the distance between \emph{normalised proper} points  \(\tb{P}\) and \(\tb{Q}\) connected by a \emph{proper} line
thanks to the relations \(|\tb{P}\cdot\tb{Q}|^2+\norm{\tb{P}\vee\tb{Q}}^2=1\) if the distance measure is elliptic
and \(|\tb{P}\cdot\tb{Q}|^2-\norm{\tb{P}\vee\tb{Q}}^2=1\) if the distance measure is hyperbolic.
Note that \(\norm{\tb{a}\cdot\tb{P}}=1\) and \(|\tb{P}\cdot\tb{Q}|=1\) in degenerate spaces and therefore
the inner product cannot be used to compute the distance there.

The metric structure of non-kinematic spaces is relatively simple, since 
the distance from a given proper point is defined to any other proper point.
In elliptic space \El{2}, all points are proper (including points represented by stacks) and all lines are proper, 
so the distance is defined for any two points in the space
and the angle is defined for any two lines.
In Euclidean space \E{2}, points represented by sheaves are proper, points represented by stacks are null,
and the only null line is \(\e_0\). 
Even though the distance to null points is not defined, one can say informally that the null points are at infinity
and \(\e_0\) is the line at infinity.
In hyperbolic space \Hy{2}, null points lie on the unit circle defined by \(x^2+y^2=1\), proper points are inside the unit circle,
and improper points are outside of it.
Proper lines pass through proper points, null lines are tangent to the unit circle, and all other lines are improper.
The unit circle represents infinity in hyperbolic space and the distance is defined only for the points  \((x,y)\) which satisfy \(x^2+y^2<1\).
The null points on a proper line \(\tb{a}\) in \Hy{2} lie on the unit circle and
are given by \((\tb{a}\pm\I)\tb{b}\), where \(\tb{a}\) is assumed to be normalised and 
\(\tb{b}\) is any non-null line passing through \(\tb{a}\I\).

\begin{figure}[t!]
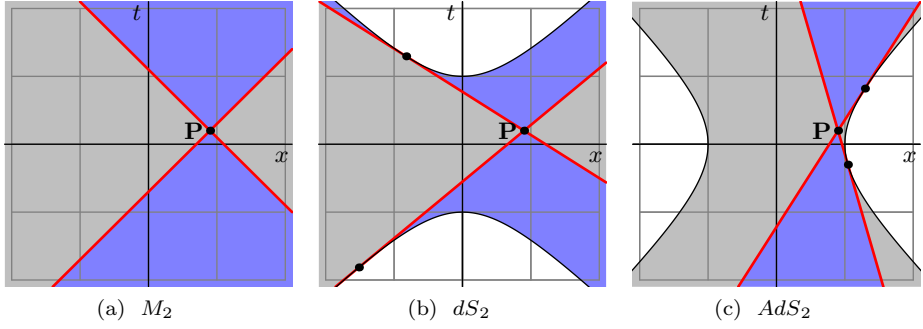

\begin{subfloatenv}{ \M{2}}
\begin{asy}
import Figure2D;
metric = Metric(Minkowski);
real w = 2.0, m = 0.1, unit_size=0.9;
Figure f = Figure(w,m,xaxis_name="$x$",yaxis_name="$t$", unit_size=unit_size);

MV P = Point(1,1/8,3/4);
MV P = Point(1,9/10,2/10);
P/=norm(P);
MV n1 = join(P, (P-1)*cross(P,e_1*e_2));
MV n2 = join(P, (P+1)*cross(P,e_1*e_2));

MV a1 = Line(w+m,1,0);
MV Q1 = wedge(a1,n1);
MV R1 = wedge(a1,n2);

MV a2 = Line(2*(w+m),-1,0);
MV R2 = wedge(a2,n1);
MV Q2 = wedge(a2,n2);

pair p(MV P) {return topair(P/P.w); }

fill(p(P)--p(Q1)--p(R1)--cycle,mediumgray);
fill(p(P)--p(R2)--p(Q2)--cycle,mediumgray);
fill(p(R2)--p(P)--p(R1)--cycle,lightblue);
fill(p(Q2)--p(P)--p(Q1)--cycle,lightblue);
Figure(w,m,xaxis_name="",yaxis_name="", unit_size=unit_size);

f.line(n1,draw_orientation=false,pen=red+1);
f.line(n2,draw_orientation=false,pen=red+1);
f.point(P,"$\textbf{P}$",align=(-0.5,0),draw_orientation=false);

f.crop();

\end{asy}
\end{subfloatenv}%
\begin{subfloatenv}{ \dS{2}}
\begin{asy}
import Figure2D;
metric = Metric(deSitter);
real w = 2.0, m = 0.1, unit_size=0.9;
Figure f = Figure(w,m,xaxis_name="$x$",yaxis_name="$t$", unit_size=unit_size);

MV P = Point(1,1/8,3/4);
MV P = Point(1,9/10,2/10);
P/=norm(P);
MV N1 = (1-P)*cross(e_1*e_2,P);
MV N2 = (1+P)*cross(e_1*e_2,P);
MV n1 = join(P, (P-1)*cross(P,e_1*e_2));
MV n2 = join(P, (P+1)*cross(P,e_1*e_2));

void draw_f(real f(real), real x1, real x2=-x1) { draw(graph(f,x1,x2,operator ..)); }
real f(real x) { return sqrt(1+x^2); }
real g(real x) { return -f(x); }

MV a1 = Line(w+m,1,0);
MV Q1 = wedge(a1,n1);
MV R1 = wedge(a1,n2);

MV a2 = Line(w+m,-1,0);
MV R2 = wedge(a2,n1);
MV Q2 = wedge(a2,n2);

pair p(MV P) {return topair(P/P.w); }

path upper = graph(f,-(w+m),w+m,operator ..);
path lower = graph(g,-(w+m),w+m,operator ..);

fill(p(P)--p(Q1)--p(R1)--cycle,mediumgray);
fill(p(P)--p(R2)--p(Q2)--cycle,mediumgray);
fill(lower--p(R2)--p(P)--p(R1)--cycle,lightblue);
fill(upper--p(Q2)--p(P)--p(Q1)--cycle,lightblue);
Figure(w,m,xaxis_name="",yaxis_name="", unit_size=unit_size);

f.line(n1,draw_orientation=false,pen=red+1);
f.line(n2,draw_orientation=false,pen=red+1);
f.point(P,"$\textbf{P}$",align=(-0.5,0),draw_orientation=false);
f.point(N1," ",draw_orientation=false);
f.point(N2," ",draw_orientation=false);

draw_f(f, w);
draw_f(g, w);

f.crop();

\end{asy}
\end{subfloatenv}%
\begin{subfloatenv}{ \AdS{2}}
\begin{asy}
import Figure2D;
metric = Metric(AntideSitter);
real w = 2.0, m = 0.1, unit_size=0.9;
Figure f = Figure(w,m,xaxis_name="$x$",yaxis_name="$t$", unit_size=unit_size);

MV P = Point(1,1/8,3/4);
MV P = Point(1,9/10,2/10);
P/=norm(P);
MV N1 = (1-P)*cross(e_1*e_2,P);
MV N2 = (1+P)*cross(e_1*e_2,P);
MV n1 = join(P, (P-1)*cross(P,e_1*e_2));
MV n2 = join(P, (P+1)*cross(P,e_1*e_2));

void draw_f(real f(real), real x1, real x2=-x1) { draw(reflect((-1,-1),(1,1))*graph(f,x1,x2,operator ..)); }
real f(real x) { return sqrt(1+x^2); }
real g(real x) { return -f(x); }

MV a1 = Line(w+m,0,1);
MV Q1 = wedge(a1,n1);
MV R1 = wedge(a1,n2);

MV a2 = Line(w+m,0,-1);
MV R2 = wedge(a2,n1);
MV Q2 = wedge(a2,n2);

pair p(MV P) {return topair(P/P.w); }

path upper = reflect((-1,-1),(1,1))*graph(f,-(w+m),w+m,operator ..);
path lower = reflect((-1,-1),(1,1))*graph(g,-(w+m),w+m,operator ..);

fill(p(P)--p(Q1)--p(R1)--cycle,lightblue);
fill(p(P)--p(R2)--p(Q2)--cycle,lightblue);
fill(lower--p(R2)--p(P)--p(R1)--cycle,mediumgray);
fill(upper--p(Q2)--p(P)--p(Q1)--cycle,mediumgray);
Figure(w,m,xaxis_name="",yaxis_name="", unit_size=unit_size);

f.line(n1,draw_orientation=false,pen=red+1);
f.line(n2,draw_orientation=false,pen=red+1);
f.point(P,"$\textbf{P}$",align=(-0.5,0),draw_orientation=false);
f.point(N1," ",draw_orientation=false);
f.point(N2," ",draw_orientation=false);

draw_f(f, w);
draw_f(g, w);
f.crop();

\end{asy}
\end{subfloatenv}
\caption{The extent of the kinematic space consisting of those points
to which the distance from \(\tb{P}=\e_{12}+\tfrac{9}{10}\e_{20}+\tfrac{2}{10}\e_{01}\) is defined (blue areas); 
other proper points (grey areas), null points (black curves), 
and null lines passing through \(\tb{P}\) (red lines).}
\label{metric content of kinematic spaces}
\end{figure}

In kinematic spaces, the situation is more complicated.
To emphasize the distinction between kinematic and non-kinematic spaces, 
I will replace \(y\) with \(t\) and \(b\) with \(h\), so that a point and a line 
can be written as \(\tb{P}=w\e_{12}+x\e_{20}+t\e_{01}\) and \(\tb{a}=d\e_0+a\e_1+h\e_2\).
In Minkowski space, points represented by sheaves are proper, points represented by stacks are null,
and the line \(\e_0\) is also null as in Euclidean space.
However, Minkowski space possesses many other null lines besides \(\e_0\).
The diagonal lines \(\tb{a}=d\e_0+a\e_1\pm a\e_2\), with  \(|a|=|h|\), are all null.
A line is proper if its slope is steeper (\(|a|>|h|\)) than that of the diagonals and improper if its slope is shallower (\(|a|<|h|\)).
The extent of the space, consisting of those points to which the distance from a given proper point \(\tb{P}\) is defined,
is bound by the two null lines passing through \(\tb{P}\) (see Figure~\ref{metric content of kinematic spaces}(a)).
The distance to the null points is undefined, but one can say informally that a null point \(\tb{Q}\)
is at the infinite distance from \(\tb{P}\) provided that \(\tb{P}\vee\tb{Q}\) has a steeper slope than the null lines passing through  \(\tb{P}\).

In de-Sitter space, a point \(\tb{P}=\e_{12}+x\e_{20}+t\e_{01}\)
is null if it lies on a hyperbola defined by \(t^2-x^2=1\), which I will call the null hyperbola in \dS{2}.
The point is proper if \(t^2-x^2<1\) and improper if \(t^2-x^2>1\).
In addition to the null points represented by sheaves, 
there are two null points represented by stacks and given by \(\tb{P}=\e_0\wedge(\e_1\pm \e_2)\).
One of these null points lies on the line \(\e_1+\e_2\) and the other lies on \(\e_1-\e_2\),
which are the two diagonal lines passing through the origin.
A point \(\tb{P}=\e_0\wedge(a\e_1+h\e_2)\) is proper if  \(|a|<|h|\) and improper if  \(|a|>|h|\),
e.g.\ the point \(\e_{20}\), which lies on the \(x\)-axis, is proper.
There are two null lines passing through any proper point \(\tb{P}\).
They are both tangent to the null hyperbola and are given by \((\tb{P}\pm1)\tb{b}\), where \(\tb{P}\) is assumed to be normalised
and \(\tb{b}\) is any non-null line passing through \(\tb{P}\).
For instance, the null lines passing through the origin are given by  \(\e_1+\e_2\) and \(\e_1-\e_2\).
Likewise, the null points where the tangent lines touch the null hyperbola are given by \((\tb{P}\pm1)\tb{Q}\),
where \(\tb{Q}\) is any non-null point on the line \(\tb{P}\I\).
Proper lines passing through \(\tb{P}\) are bound by the null lines through \(\tb{P}\) and
the extent of the space, consisting of those points to which the distance from a given proper point \(\tb{P}\) is defined,
is bound by the null lines through \(\tb{P}\) and the null hyperbola (see Figure~\ref{metric content of kinematic spaces}(b)).
Note that \(\e_0\) is a proper line in \dS{2}.
The null points on a proper line \(\tb{a}\) lie on the null hyperbola and
are given by the same formulas as in hyperbolic space.

Anti de-Sitter space has a structure somewhat similar to the structure of de-Sitter space (see Figure~\ref{metric content of kinematic spaces}(c)).
In anti de-Sitter space, a point \(\tb{P}=\e_{12}+x\e_{20}+t\e_{01}\)
is null if it lies on a hyperbola defined by \(t^2-x^2=-1\), which may be called the null hyperbola in \AdS{2}.
The point is proper if \(t^2-x^2>-1\) and improper if \(t^2-x^2<-1\).
The two points \(\tb{P}=\e_0\wedge(\e_1\pm \e_2)\) are null, as in de-Sitter space.
A point \(\tb{P}=\e_0\wedge(a\e_1+h\e_2)\) is proper if  \(|a|>|h|\) and improper if  \(|a|<|h|\),
e.g.\ the point \(\e_{01}\), which lies on the \(t\)-axis, is proper.
The null lines passing through a proper point \(\tb{P}\) are tangent to the null hyperbola and are given by the same formulas
as in de-Sitter space.
Proper lines passing through \(\tb{P}\) are bound by the null lines through \(\tb{P}\) in the same fashion as in de-Sitter space.
Note that \(\e_0\) is an improper line in \AdS{2}.
Proper lines do not intersect the null hyperbola in \AdS{2}.

It is instructive to derive explicit expressions for the distance from the origin \((0,0)\) to another proper point.
In non-kinematic spaces, I get the following expressions:
\begin{equation}
r=\sqrt{x^2+y^2}, \quad
\sin r = \frac{\sqrt{x^2+y^2}}{\sqrt{1+x^2+y^2}}, \quad
\sinh r = \frac{\sqrt{x^2+y^2}}{\sqrt{1-(x^2+y^2)}},
\end{equation}
for the distance \(r\) from the origin to \((x,y)\) in \E{2}, \El{2}, \Hy{2}, respectively.
And in kinematic spaces, I get
\begin{equation}
r=\sqrt{t^2-x^2}, \quad
\sin r = \frac{\sqrt{t^2-x^2}}{\sqrt{1+(t^2-x^2)}}, \quad
\sinh r = \frac{\sqrt{t^2-x^2}}{\sqrt{1-(t^2-x^2)}}
\end{equation}
for the distance from the origin to \((x,t)\) in \M{2}, \AdS{2}, \dS{2}, respectively.

The projective model allows one to visualise points of the two-dimensional space represented by sheaves as a plane.
The drawback of this visualisation is that points represented by stacks cannot be shown directly.
The advantage is that only two dimensions are required for visualisation and all geodesics are represented by straight lines.
It is common to visualise proper points of a non-degenerate space as a suitable curved surface embedded in \R{3}.
One can define the surface as a set of coordinates of all normalised proper points.
In elliptic and hyperbolic spaces, it consists of all triples \((w,x,y)\) which satisfy \(\tb{P}^2=-1\), where \(\tb{P}=w\e_{12}+x\e_{20}+y\e_{01}\).
This yields a unit sphere, \(w^2+x^2+y^2=1\), for elliptic space, where all points are proper,
and a hyperboloid of two sheets, \(w^2-(x^2+y^2)=1\), for the proper points of hyperbolic space.
In de-Sitter and anti de-Sitter spaces, the surface consists of all triples \((w,x,t)\) which satisfy \(\tb{P}^2=1\), 
 where \(\tb{P}=w\e_{12}+x\e_{20}+t\e_{01}\).
This yields a hyperboloid of one sheet, \(w^2-(t^2-x^2)=1\), for the proper points of de-Sitter space and
another hyperboloid of one sheet, \(w^2+(t^2-x^2)=1\), for the proper points of anti de-Sitter space.
The advantage of such visualisation is that all proper points, including those represented by stacks in the projective model, can be shown directly.
The disadvantage is that three dimensions are required for visualisation,
improper points are not visualised at all, and geodesics can no longer be represented by straight lines.
Moreover, integration is usually required to determine distances.

\subsection{Projections and Rejections}
\label{section.projection and rejection}
In this and the next section, I consider a range of geometric transformations, which are available in the projective model thanks to the Clifford algebra.
Projections and rejections arise naturally as different components of the geometric product.
Reflections, rotations, and translations arise as a result of conjugation of lines and points by certain multivectors.
Reflections, rotations, and translations preserve distances between points and angles between lines and therefore
constitute isometries of the metric space.

The geometric product of two lines \(\tb{a}\) and \(\tb{b}\) consists of a scalar and a point, \(\tb{a}\tb{b}=\tb{a}\cdot\tb{b}+\tb{a}\wedge\tb{b}\).
If \(\tb{b}\) is invertible, then \(\tb{a}\) can be written as the sum of two lines, 
\(\tb{a}=(\tb{a}\cdot\tb{b})\tb{b}^{-1}+(\tb{a}\wedge\tb{b})\tb{b}^{-1}\).
The line \((\tb{a}\cdot\tb{b})\tb{b}^{-1}\) coincides with \(\tb{b}\), but possibly with a different weight and orientation,
and is called the projection of \(\tb{a}\) on \(\tb{b}\).
On the other hand, the line \((\tb{a}\wedge\tb{b})\tb{b}^{-1}\) passes through the point \(\tb{a}\wedge\tb{b}\) 
and the polar point \(\tb{b}\I\) of \(\tb{b}\);
it is perpendicular to \(\tb{b}\) since it passes through \(\tb{b}\I\).
It is called the rejection of \(\tb{a}\) by \(\tb{b}\).
So, any line \(\tb{a}\) can be written as the sum of its projection on and rejection by an invertible line \(\tb{b}\).
The same reasoning can be applied to the geometric product of a line and a point and the geometric product of two points.
This gives \((\tb{P}\cdot\tb{a})\tb{a}^{-1}\) for the projection of a point \(\tb{P}\) on an invertible line \(\tb{a}\)
and  \((\tb{P}\wedge\tb{a})\tb{a}^{-1}\) for the rejection of \(\tb{P}\) by \(\tb{a}\);
the projection lies on \(\tb{a}\) and the rejection takes \(\tb{P}\) to the polar point \(\tb{a}\I\) of \(\tb{a}\).
On the other hand, \((\tb{a}\cdot\tb{P})\tb{P}^{-1}\) gives the projection of \(\tb{a}\) on an invertible point \(\tb{P}\)
and \((\tb{a}\wedge\tb{P})\tb{P}^{-1}\)  gives the rejection of \(\tb{a}\) by \(\tb{P}\);
the projection passes through the points \(\tb{P}\) and \(\tb{a}\wedge(\tb{P}\I)\) and the rejection takes \(\tb{a}\) to the line \(\tb{P}\I\).
The projection \((\tb{P}\cdot\tb{Q})\tb{Q}^{-1}\) of \(\tb{P}\) on an invertible point \(\tb{Q}\) takes \(\tb{P}\) to \(\tb{Q}\)
 and the rejection \((\tb{P}\times\tb{Q})\tb{Q}^{-1}\) of \(\tb{P}\) by \(\tb{Q}\)
lies at the intersection of the lines \(\tb{P}\vee\tb{Q}\) and  \(\tb{Q}\I\).

In Euclidean space \E{2}, the rejection \((\tb{P}\wedge\tb{a})\tb{a}^{-1}\) of a finite point \(\tb{P}\) by a finite line \(\tb{a}\)
is a point at infinity in the direction perpendicular to \(\tb{a}\).
Since the rejection is at the infinite distance from \(\tb{a}\), one can combine projection and rejection with a scalar factor \(\gamma\)
to define \(\tb{P}_{\!\gamma}=(\tb{P}\cdot\tb{a})\tb{a}^{-1}+\gamma (\tb{P}\wedge\tb{a})\tb{a}^{-1}\),
which is a point obtained from \(\tb{P}\) by shifting it away (\(\gamma>1\)) or towards (\(0<\gamma<1\)) the line \(\tb{a}\)
along the direction perpendicular to \(\tb{a}\).
Moreover, if \(\tb{P}\) is at the distance \(r\) from \(\tb{a}\), then \(\tb{P}_{\!\gamma}\) is at the distance \(\gamma r\) from \(\tb{a}\).

If a point \(\tb{P}\) is invertible, then the conjugation \(\tb{P}\tb{Q}\tb{P}^{-1}\) of \(\tb{Q}\) by \(\tb{P}\) is defined.
It is a point which lies on the line \(\tb{P}\vee\tb{Q}\) and is proper if \(\tb{Q}\) is proper.
If the distance from \(\tb{P}\) to \(\tb{Q}\) is defined, then so is the distance from \(\tb{P}\) to \(\tb{P}\tb{Q}\tb{P}^{-1}\)
and they are equal.
The point \(\tb{P}\tb{Q}\tb{P}^{-1}\) is called the reflection of \(\tb{Q}\) in \(\tb{P}\).
The reflection of \(\tb{P}\) in a line \(\tb{a}\) can be identified with the reflection of \(\tb{P}\) 
in the point \((\tb{P}\cdot\tb{a})\tb{a}^{-1}\), which gives \(-\tb{a}\tb{P}\tb{a}^{-1}\).
The point  \(-\tb{a}\tb{P}\tb{a}^{-1}\) lies on the line \(\tb{a}\cdot\tb{P}\) and is proper if \(\tb{P}\) is proper.
If the distance from \(\tb{a}\) to \(\tb{P}\) is defined, then so is the distance from \(\tb{a}\) to \(-\tb{a}\tb{P}\tb{a}^{-1}\)
and they are equal as expected.

The reflection of lines can be constructed from the reflection of points as follows.
To determine the reflection of a line \(\tb{a}\) in \(\tb{P}\), 
one can express it as the join \(\tb{a}=\tb{Q}\vee\tb{R}\) of any two points on \(\tb{a}\)
and then reflect both of these points in \(\tb{P}\) and compute the join of the reflected points.
Since \((\tb{P}\tb{Q}\tb{P}^{-1})\vee(\tb{P}\tb{R}\tb{P}^{-1})=\tb{P}(\tb{Q}\vee\tb{R})\tb{P}^{-1}\),
the reflection of \(\tb{a}\) in \(\tb{P}\) is given by \(\tb{P}\tb{a}\tb{P}^{-1}\).
Similarly, I find that \(\tb{a}\tb{b}\tb{a}^{-1}\) gives the reflection of \(\tb{b}\) in \(\tb{a}\).
Observe that if \(\tb{a}\) and \(\tb{b}\) are perpendicular, 
the reflected line coincides with the original line but it has the opposite orientation,
 since \(\tb{a}\cdot\tb{b}=0\) implies \(\tb{a}\tb{b}\tb{a}^{-1}=-\tb{b}\).
This is consistent with the expectation that the bottom-up orientation of a line, 
which can be identified with the vector  \((-b,a)\) for a line \(\tb{b}=d\e_0+a\e_1+b\e_2\),
 flips when it is reflected in a perpendicular line.
I will refer to the transformations described above as the bottom-up reflections, since they are built from reflections of points.

The top-down reflection of \(\tb{b}\) in \(\tb{a}\) is given by \(-\tb{a}\tb{b}\tb{a}^{-1}\).
The minus sign is necessary for consistency with the top-down orientation, 
identified with the vector \((a,b)\) for a line \(\tb{b}=d\e_0+a\e_1+b\e_2\),
since the top-down orientation of a line is expected to remain the same when it is reflected in a perpendicular line.
The top-down reflection of \(\tb{a}\) in a point \(\tb{P}\) can be identified with the top-down reflection
of \(\tb{a}\) in the line \((\tb{a}\cdot\tb{P})\tb{P}^{-1}\), which gives \(\tb{P}\tb{a}\tb{P}^{-1}\).
The top-down reflection of points is constructed from the top-down reflection of lines.
The top-down reflection of \(\tb{P}\) in \(\tb{a}\) is given by \(\tb{a}\tb{P}\tb{a}^{-1}\),
which follows from \((-\tb{a}\tb{b}\tb{a}^{-1})\wedge(-\tb{a}\tb{c}\tb{a}^{-1})=\tb{a}(\tb{b}\wedge\tb{c})\tb{a}^{-1}\)
where \(\tb{P}=\tb{b}\wedge\tb{c}\).
Similarly, I find that the top-down reflection of \(\tb{Q}\) in \(\tb{P}\) is given by \(\tb{P}\tb{Q}\tb{P}^{-1}\).
It follows that the top-down and bottom-up reflections in a point are identical,
while the top-down and bottom-up  reflections in a line are different by a sign.

\subsection{Isometries}
\label{section.isometries}
In Euclidean space \E{2}, two consecutive reflections of a point \(\tb{P}\) in \(\tb{b}\) and then in \(\tb{a}\),
where \(\tb{a}\) and \(\tb{b}\) are normalised finite lines intersecting at a finite point, 
yield a rotation of \(\tb{P}\) around the point \(\tb{a}\wedge\tb{b}\) by twice the angle between \(\tb{a}\) and \(\tb{b}\).
This transformation can be written as \(\tb{a}\tb{b}\tb{P}\tb{b}^{-1}\tb{a}^{-1}=S\tb{P}S^{-1}\), 
where \(S=\tb{a}\tb{b}\) is an even multivector called a spinor.
In general, a multivector 
may be called a spinor if it can be written as the geometric product of an even number of normalised proper vectors,
i.e.\ \(S=\tb{a}_1\tb{a}_2\dotsm\tb{a}_{2k}\), where \(\tb{a}_i^2=1\) for all \(i=1,2,\dotsc,2k\).
Every spinor \(S\) is even, satisfies \(S\reverse{S}=1\), and has the inverse given by  \(S^{-1}=\reverse{S}\), which is also a spinor.
The geometric product of any two spinors is again a spinor.
It follows that spinors form a group. It is also a smooth manifold and therefore the group of spinors is a Lie group.
Its Lie algebra is an abstract vector space consisting of all bivectors with the addition and multiplication by scalars inherited from the Clifford algebra
and the commutator \(\times\) of the Clifford algebra serving as the product of bivectors in the Lie algebra (the Lie bracket).
Indeed, the commutator of two bivectors is again a bivector; it is antisymmetric and satisfies the Jacobi identity.

Spinors are not blades and therefore they do not represent any points or lines,
but it is convenient to refer to them as geometric objects in \T{2}.
There is a close connection between spinors and points in \T{2}.
Any multivector that can be written as \(e^\tb{P}\) for some point \(\tb{P}\) is a spinor.
The reverse holds in elliptic space only, i.e.\ every spinor \(S\) can be written as \(S=e^\tb{P}\) for some point \(\tb{P}\) in \El{2}.
In the other spaces, a spinor can only be written as either \(S=e^\tb{P}\) or \(S=-e^\tb{P}\) for some point \(\tb{P}\) 
(see \cite[Section 5.3]{gunn2011geometry} and \cite[Sections 3-5, 3-8]{hestenes1987clifford}).
For example, \(S=(\e_0-\e_2)\e_2=-e^{\e_{20}}\) is a spinor in \E{2} but \(S\ne e^{\tb{P}}\) for any point \(\tb{P}\).
The action of \(S\) on a proper point or line via conjugation yields a transformation that preserves distances and angles 
(see \ref{appendix.isometries}).
The conjugation by \(S=e^\tb{P}\) yields a rotation around \(\tb{P}\) if it is proper and a translation if it is null or improper.
Since there are no null or improper points in \El{2}, no  translations are possible there.
Note that the action of \(-e^\tb{P}\) is equivalent to that of \(e^\tb{P}\), so there is no need to consider spinors \(-e^\tb{P}\).

In non-kinematic spaces, the terms \emph{rotation} and \emph{translation} have the expected meaning.
In  Euclidean space \E{2}, the conjugation \(S\tb{P}S^{-1}\) of a finite point \(\tb{P}\) by \(S=e^{-\frac{1}{2}\alpha\tb{R}}\), 
where \(\alpha>0\) is a scalar and \(\tb{R}=\e_{12}+x\e_{20}+y\e_{01}\), 
yields a counterclockwise rotation of \(\tb{P}\) around \(\tb{R}\) by the angle \(\alpha\);
a clockwise rotation obtains if \(\alpha<0\).
The conjugation \(T\tb{P}T^{-1}\) of a finite point \(\tb{P}\) by \(T=e^{-\frac{1}{2}r\tb{T}}\),
where \(r>0\) is a scalar and \(\tb{T}=\e_0\wedge(\e_1\cos\alpha+\e_2\sin\alpha)\) is a point at infinity, 
yields a translation of \(\tb{P}\) in the direction \((\cos\alpha,\sin\alpha)\) by the distance \(r\).
The same spinors can be applied to a finite line, which yields the same transformations, i.e.\ 
every point on the line experiences rotation or translation as described above.
Rotations and translations  constitute the group of rigid body motions of the two-dimensional Euclidean space
and every rigid body motion can be represented by the action of a suitable spinor.

In Minkowski space \M{2}, a rotation around a proper point 
may be interpreted in kinematic terms as the addition of velocities (or boost)
when applied to proper lines or the Lorentz transformation when applied to proper points.
Rotations together with translations in \M{2} constitute the proper Poincar\'e group.
For instance, a boost (along the \(x\)-axis) by the speed \(v=\tanh\phi_b\)
is given by the conjugation \(B\tb{a}B^{-1}\), where \(B=e^{\frac{1}{2}\phi_b\e_{12}}\) and
\(\tb{a}=\e_1\cosh\phi-\e_2\sinh\phi\) is a proper line, which can be identified with a
world-line of a point moving along the \(x\)-axis at the speed \(u=\tanh\phi\).
The parameter \(\phi\in\R{}\) is called rapidity and \(|\phi|\) gives the angle between \(\tb{a}\) and \(\e_1\), the \(t\)-axis.
The Lorentz transformation consistent with the boost is given by \(B\tb{P}\!B^{-1}\), 
where \(\tb{P}=\e_{12}+x\e_{20}+t\e_{01}\) is a proper point identified with the event \((x,t)\).
On the other hand, the conjugations \(N\tb{P}N^{-1}\) of a point \(\tb{P}\) by a spinor \(N=e^{-\frac{1}{2}\lambda\tb{N}}\), 
where \(\lambda\in\R{}\) and \(\tb{N}\) is a null point, yields a translation, e.g.\
a translation by \(\lambda>0\) along the \(t\)-axis in the direction of increasing \(t\) 
is induced by the spinor \(N=e^{-\frac{1}{2}\lambda\e_{20}}\).

\section{Discussion}
The model described in this paper pertains to the two-dimensional spaces, but it can be readily extended to higher dimensions.
The construction is based on projective duality and follows the recipe used in two dimensions.
The notion of projective duality is modified according to the dimension of the space.
For instance, in a three-dimensional space the defining equation reads \(1+ax+by+cz=0\) and, therefore, a point is dual to a plane.
Planes are at the top of the hierarchy and the top-down model is based on planes.
A finite line is understood as a sheaf of all planes passing through the line and a point is a bundle of all planes passing through the point.
If the space is Euclidean, stacks of parallel planes represent lines at infinity, 
while points at infinity can be thought of as stacks of parallel lines.
The target space is embedded into a four-dimensional model space \R{4} 
and geometric objects are identified with linear subspaces of \R{4} and its dual \R{4*},
e.g.\ planes, lines, and points are respectively represented by vectors, bivectors, and trivectors of \(\bigwedge\R{4*}\).
The structure of \(\bigwedge\R{4*}\) is more complex than \(\bigwedge\R{3*}\), 
which reflects a more complex structure of the three-dimensional space.
The Grassmann algebra is 16-dimensional as an abstract vector space and contains non-simple bivectors,
which correspond to interesting geometric structures in the target space.

As the dimension of the target space increases, the structure of the Clifford algebra used to describe it adjusts to account for 
the geometric structures possible in a higher-dimensional environment.
The geometric tool-set provided by a higher-dimensional Clifford algebra is correspondingly richer than in lower dimensions,
but it has the same basic characteristics.
Projections, rejections, and reflections are constructed in higher dimensions by the same recipes.
Spinors are defined in the same way via an even number of reflections in the geometric objects at the top of the hierarchy,
e.g.\ planes in \T{3} and hyperplanes in \T{4}.
The group of rigid body motions is richer as it includes other transformations besides the usual rotations and translations.
Furthermore, a rotation no longer occurs around a point as in \T{2}  but around a line in \T{3} and around a plane in \T{4},
which is reflected in the fact that bivectors represent lines in \T{3} and planes in \T{4}.
The extension of the projective model described in this paper to three- and four-dimensional spaces will be 
given in a separate paper (three-dimensional non-kinematic spaces are described in \cite{gunn2011geometry}).

It is possible to transfer the metric from the dual to the model space in a way consistent with the Cayley-Klein construction \cite{gunn2011geometry}.
For instance, the inner product of two points \(\tb{P}=w_p\e_{12}+x_p\e_{20}+y_p\e_{01}\) and \(\tb{Q}=w_q\e_{12}+x_q\e_{20}+y_q\e_{01}\) in \Hy{2}
is given by \(\tb{P}\cdot\tb{Q}=-w_pw_q+x_px_q+y_py_q\) and hence one can define
\(\tb{x}\cdot\tb{y}=-w_pw_q+x_px_q+y_py_q\) for \(\tb{x}=w_p\e^{0}+x_p\e^{1}+y_p\e^{2}\) and \(\tb{y}=w_q\e^{0}+x_q\e^{1}+y_q\e^{2}\),
which yields a Clifford algebra over \R{3} with the properties identical to those of the Clifford algebra over \R{3*}.
However, the transferred metric and the resulting Clifford algebra retain their utility only if the original metric is non-degenerate.
In the case of degenerate spaces such as Euclidean and Minkowski spaces, 
the transferred metric is too degenerate and is inconsistent with the expected metric properties of these spaces.
Nevertheless, there are homogeneous models (e.g.\ \cite[Chapter 11]{dorst2009geometric}) that apply in a limited way
a Clifford algebra over \R{n+1} (with a non-degenerate metric) to the \(n\)-dimensional Euclidean geometry, since
some geometric transformations such as rotations around the origin yield equivalent results in all non-kinematic spaces
and therefore can be used to represent Euclidean rotations around the origin.

The conformal model introduced in \cite{hestenes2001new} (see also \cite{dorst2009geometric}) 
accomplishes some of the same objectives as the projective model in providing an algebraic framework for
representing geometric objects in Euclidean space and their transformations.
Unlike the projective model, it requires two extra dimensions for the model space, 
so that a Clifford algebra over \R{5} (with a non-degenerate metric) is needed to represent geometry in \E{3}.
The payback is that in addition to representing points, lines, and planes, it can be used to represent circles and spheres with blades.
It also provides a convenient way of representing scaling via conjugation by certain spinors.
The comparative advantages of the projective model are computational efficiency and conceptual simplicity.
As pointed out in \cite{gunn2011geometry} the projective model provides a perfect framework for rigid body dynamics, 
whereas the extra degrees of freedom inherent in the conformal model present a liability and may need to be eliminated \cite{lasenby2011rigid}.
Moreover, the projective model provides the simplest representation of the Poincar\'e group and its Lie algebra,
which has a significant potential of further applications in the study of relativistic space-times and other theories that rely on such space-times.

One can perform calculations in Clifford algebra by hand, but it may be onerous and error-prone.
Fortunately, there is a growing range of software based on Grassmann and Clifford algebras:
for symbolic algebra \cite{browne_grassmann,ablamowitz_CLIFFORD,bromborsky_sympy},
visualisation \cite{perwass_CLUCalc,dorst_GAViewer}, and computing in C++ and Java
(see \cite{hildenbrand2013foundations} for further details).
Various implementation issues are also discussed in \cite{dorst2009geometric}.

\section{Conclusion}
Projective geometry and Clifford algebra give rise to the projective model of homogeneous metric spaces known as Cayley-Klein geometries.
The projective model provides a set of tools for manipulating geometric objects such as points, lines, and planes.
These tools are based on a Clifford algebra, which replaces synthetic reasoning and trigonometry with geometric  multiplication.
The projective model obviates the need for frequent coordinate transformations, since
one can manipulate geometric objects directly, while the coordinates are dealt with automatically by the geometric product.
Moreover, Clifford algebra takes into account all metric aspects automatically thereby freeing the user to 
apply the same tool-set in any metric space.

Synthetic techniques are well developed for Euclidean geometry in two and three dimensions and for spherical geometry,
but they do not extend easily to higher dimensions and to other homogeneous metric spaces.
The advantage of the projective model is that it shields the user from the intricate and intractable details of the underlying geometric structures,
which one has to deal with explicitly in the synthetic approach.
The projective model gives direct access to various geometric objects and essentially turns geometry into algebra
by replacing geometric manipulations with algebraic operations, 
enabling one to efficiently solve geometric problems in higher dimensions and in spaces with a non-Euclidean metric;
in particular, the study of relativistic kinematic spaces is simplified.
Moreover, it offers a novel way of describing the Poincar\'e group and its Lie algebra, which is simpler than the traditional matrix methods
as it does not rely on coordinates.

\ack
The author thanks Charles Gunn for clarifying certain aspects of the projective model.

\appendix 
\section{Duality transformation}
\label{appendix.duality Rn}

For a basis multivector \(\e_{I}\), where \(I\) is a list of indexes for this basis multivector, I define \(\J:\bigwedge\R{(n+1)*}\to\bigvee\R{n+1}\) by
\begin{equation} 
\J(\e_I)=\e^{\reverse{I^\perp}}\label{J},
\end{equation}
where \({I}^\perp\) denotes a list of indexes complementary to \(I\), such that the concatenation \({I}^\perp I\) is an even permutation of \(01\dots n\),
and \(\reverse{I^{\perp}}\) denotes the reversion of the indexes in \(I^{\perp}\),  e.g.\ \(\reverse{012}=210\)
(it is assumed that \(\e_{\varnothing}=1\), \(\e^{\varnothing}=1\), \(\reverse{\varnothing}= \varnothing\), and 
\(\reverse{k}=k\) for a list consisting of a single integer \(k\)).
This definition extends to arbitrary multivectors in \(\bigwedge\R{(n+1)*}\) by linearity.
Applying (\ref{J}) when \(n=2\) gives Table~\ref{duality R3}.

If the metric is elliptic, then
\begin{equation}
\J(A)=\Id(A\I^{-1})
\label{duality elliptic}
\end{equation} for any multivector \(A\in\bigwedge\R{(n+1)*}\),
where \(\Id:\bigwedge\R{(n+1)*}\to\bigvee\R{n+1}\) is defined on the basis multivectors by \(\Id(\e_I)=\e^I\) and extended to 
general multivectors by linearity (it is called Riesz map in \cite{roman2008advanced}).
Indeed, since \({I}^\perp I\) is an even permutation of \(01\dots n\),  I can write
\(\I=\e_{01\dots n}=\e_{I^{\perp}I}\) for any basis blade \(\e_I\) and therefore 
\begin{equation}
\notag
\e_I\I^{-1}=
\e_I\reverse{\I}=
\e_I\reverse{\e}_{I^{\perp}I}=
\e_I\reverse{\e}_I \reverse{\e}_{I^{\perp}}=
\reverse{\e}_{I^{\perp}}=
\e_{\reverse{I^{\perp}}},
\end{equation}
which gives \(\J(\e_I)=\Id(\e_I\I^{-1})\), and (\ref{duality elliptic}) follows from linearity of \(\J\) and \(\Id\).
The duality transformation does not depend on a metric, but one can use (\ref{duality elliptic}) in calculations as long as the elliptic metric is used.
Hence, the inverse transformation \(\J^{-1}:\bigvee\R{n+1}\to\bigwedge\R{(n+1)*}\) can be written as
\(\J^{-1}(\e^I)=\Id^{-1}(\e^I)\I=\e_I \I^{-1}\I^{2}\),
and since \(\I^2=-1\) if \(n=1,2\) and \(+1\) if \(n=3,4\) when the metric is elliptic, it follows that
\begin{equation}
\J^{-1}(\e^I)
=
\left\{
\begin{aligned}
& -\e_{\reverse{I^\perp}},\quad\textrm{if }n=1,2\\
& +\e_{\reverse{I^\perp}},\quad\textrm{if }n=3,4\\
\end{aligned}
\right.
\end{equation}
This can be used to streamline calculations of the inverse.
If \(n=3\) or \(4\), then \(\J^{-1}\) can be obtained from \(\J\) by lowering and raising the indexes, 
e.g.\ \(\J(\e_{230})=\e^{41}\) gives \(\J^{-1}(\e^{230})=\e_{41}\) for \(n=4\).
If \(n=1\) or \(2\), then in addition to lowering and raising the indexes it is also necessary to change the sign,
e.g.\ \(\J(\e_{12})=\e^0\) gives \(\J^{-1}(\e^{12})=-\e_0\) for \(n=2\).
The definition~(\ref{J}) is not sufficient if \(n=1\), because 01 is the only possible even permutation of 01, which gives 
\(\J(1)=-\e^{01}\),
\(\J(\e_1)=\e^0\),
\(\J(\e_{01})=1\).
I can then appeal to (\ref{duality elliptic}) to define  \(\J(\e_0)=-\e^1\).

Note that \cite{gunn2011geometry} uses Poincar\'e duality \(\mathcal{P}\) \cite{greub1967multilinear}, 
which  satisfies \(\mathcal{P}(\e_I)=\e^{I^\perp}\) on a basis multivector \(\e_I\),
to define the meet and join. 
In general,  \(\mathcal{P}(A)\ne\J(A)\) for a  multivector \(A\in\bigwedge\R{(n+1)*}\),
but if \(A\) is a simple \(k\)-vector, then \(\mathcal{P}(A)=\pm\J(A)\), where the sign depends on \(k\) and \(n\).
The transformation \(\mathcal{P}\) does not satisfy (\ref{duality elliptic}),
but \(\mathcal{P}\) is related to Hodge duality (also called Hodge star) by
\(\mathcal{P}(A)=\Id(\mathcal{H}(A))\), where \(\mathcal{H}:\bigwedge\R{(n+1)*}\to\bigwedge\R{(n+1)*}\) 
denotes the Hodge duality transformation based on the elliptic metric.

\section{Isometries}
\label{appendix.isometries}
A linear transformation \(\ts{f}:\R{3*}\to\R{3*}\) is defined on vectors of \R{3*}, or lines in \T{2},
and satisfies \(\ts{f}(\tb{a}+\tb{b})=\ts{f}(\tb{a})+\ts{f}(\tb{b})\) and \(\ts{f}(s\tb{a})=s\ts{f}(\tb{a})\).
It can be extended to general multivectors of \(\bigwedge\R{3*}\) by means of 
\begin{equation}
\ts{f}(s)=s, \quad
\ts{f}(\tb{a}\wedge\tb{b})=\ts{f}(\tb{a})\wedge\ts{f}(\tb{b}), \quad
\ts{f}(\tb{a}\wedge\tb{b}\wedge\tb{c})=\ts{f}(\tb{a})\wedge\ts{f}(\tb{b})\wedge\ts{f}(\tb{c}),
\label{f on blades}
\end{equation} 
where \(s\) is a scalar and \(\tb{a}, \tb{b}, \tb{c}\) are vectors in \R{3*},
which implies \(\ts{f}(A\wedge B)=\ts{f}(A)\wedge\ts{f}(B)\) for any \(A,B\in\bigwedge\R{3*}\).
The determinant of \(\ts{f}\) satisfies \(\ts{f}(\I)=\det{(\ts{f})}\I\), where \(\I=\e_{012}\) is the basis pseudoscalar of \(\bigwedge\R{3*}\).
The adjoint \(\bar{\ts{f}}:\R{3}\to\R{3}\) is defined implicitly by \(\ts{f}(\tb{a})[\tb{x}]=\tb{a}[\bar{\ts{f}}(\tb{x})]\),
where \(\tb{a}\) and \(\ts{f}(\tb{a})\) are treated as linear functionals,
and is characterised by the property that its matrix representation in the standard basis equals the transpose of the matrix representation of \(\ts{f}\).
A linear transformation \(\ts{g}:\R{3}\to\R{3}\) satisfies the same properties and can be extended to multivectors of \(\bigvee\R{3}\)
by using the outer product \(\vee\) in the same fashion  \(\wedge\) is used in (\ref{f on blades}):
\begin{equation}
\ts{g}(s)=s,\quad
\ts{g}(\tb{x}\vee\tb{y})=\ts{g}(\tb{x})\vee\ts{g}(\tb{y}),\quad
\ts{g}(\tb{x}\vee\tb{y}\vee\tb{z})=\ts{g}(\tb{x})\vee\ts{g}(\tb{y})\vee\ts{g}(\tb{z}),
\end{equation}
where \(s\) is a scalar and \(\tb{x},\tb{y},\tb{z}\) are vectors in \(\R{3}\).
Similarly, the equality \(\ts{g}(X\vee Y)=\ts{g}(X)\vee \ts{g}(Y)\) holds for any \(X,Y\in\bigvee\R{3}\).
If \(\ts{f}\) is invertible, its adjoint \(\bar{\ts{f}}\) is also invertible and 
\begin{equation}
\J(\ts{f}(A)) = \det(\ts{f})\bar{\ts{f}}^{-1}(\J(A))
\label{adjoint and inverse}
\end{equation}
holds for any \(A\in\bigwedge\R{3*}\) (cf.\ \cite[Section 3-1]{hestenes1987clifford}),
where \(\ts{f}\) and \(\bar{\ts{f}}^{-1}\) are assumed to be extended to general multivectors of \(\bigwedge\R{3*}\) and \(\bigvee\R{3}\), respectively.
It follows that
\begin{equation}
\ts{f}(A\vee B)=\ts{f}(A)\vee\ts{f}(B)/\det(\ts{f})
\end{equation}
for any \(A,B\in\bigwedge\R{3*}\), provided that \(\ts{f}\) is invertible.
If \(\ts{f}\) is a reflection, rotation, or translation, then \(\det(\ts{f})=\pm1\) and \(\norm{\ts{f}(A)}=\norm{A}\) for any \(A\in\bigwedge\R{3*}\)
(\(\det(\ts{f})=1\) for rigid body motions such as rotations and translations).
Therefore, I get
\begin{equation}
\notag
\norm{\tb{P}\vee\tb{Q}}=
\norm{\ts{f}(\tb{P}\vee\tb{Q})}=
\norm{\ts{f}(\tb{P})\vee\ts{f}(\tb{Q})}
\end{equation}
and
\begin{equation}
\notag
|\tb{a}\cdot\tb{b}|=
|\ts{f}(\tb{a}\cdot\tb{b})|=
|\ts{f}((\tb{a}\I)\vee\tb{b})|=
|\ts{f}(\tb{a}\I)\vee\ts{f}(\tb{b})|=
|(\ts{f}(\tb{a})\I)\vee\ts{f}(\tb{b})|=
|\ts{f}(\tb{a})\cdot\ts{f}(\tb{b})|
\end{equation}
for any points \(\tb{P}\) and \(\tb{Q}\) and lines \(\tb{a}\) and \(\tb{b}\).
It follows that reflections, rotations, and translations preserves distances and angles, i.e.\ these transformations are isometries.

\section*{References}
\bibliographystyle{elsarticle-num}
\bibliography{p.bib}

\end{document}